\documentclass[10pt,a4paper]{article}

\usepackage{amsmath,amssymb, bbm}
\usepackage{color}
\newcommand{\C}{\mathbb{C}}
\newcommand{\R}{\mathbb{R}}

\newcommand{\Z}{\mathbb{Z}}

\newcommand{\B}{\mathbb{B}}

\newcommand{\Bo}{\mathbb{B}_0}
\newcommand{\Co}{\mathbb{C}_0}
\newcommand{\Do}{\mathbb{D}_0}

\def\cA{{\mathcal A}}
\def\cC{{\mathcal C}}
\def\cB{{\mathcal B}}

\def\cP{{\mathcal P}}
\def\cQ{{\mathcal Q}}
\newcommand{\ee}{\varepsilon}
\renewcommand{\aa}{\alpha}
\newcommand{\aA}{\alpha_1}
\newcommand{\aB}{\alpha_2}
\newcommand{\aC}{\alpha_3}

\renewcommand{\div}{{\rm div}\,}

\newcommand{\Sum}{\displaystyle \sum}
\newcommand{\Hs}{\dot{H^s}}
\newcommand{\bd}{\beta \delta}

\def\d{\partial}

\def\ddj{\dot \Delta_j}

\def\ddk{\dot \Delta_k^v}
\def\tilde{\widetilde}
\def\hat{\widehat}
\newcommand{\D}{\Delta}
\newcommand{\DF}{\Delta_F}

\newcommand{\n}{\nabla}
\newcommand{\G}{\Gamma}

\newcommand{\Fe}{F_{ext}}
\newcommand{\Om}{\Omega}
\newcommand{\Ome}{\Omega_\varepsilon}
\newcommand{\tOm}{\tilde{\Omega}_{QG}}

\newcommand{\ve}{v_\ee}
\newcommand{\Ue}{U_\ee}
\newcommand{\Uqg}{U_{\ee,QG}}
\newcommand{\Uosc}{U_{\ee, osc}}
\newcommand{\Uoe}{U_{0,\ee}}
\newcommand{\Uoqg}{U_{0,\ee,QG}}
\newcommand{\Uoosc}{U_{0,\ee, osc}}
\newcommand{\tUqg}{\tilde{U}_{QG}}
\newcommand{\tvqg}{\tilde{v}_{QG}}
\newcommand{\tUoqg}{\tilde{U}_{0, QG}}

\newcommand{\We}{W_\ee}
\newcommand{\Wet}{W_\ee^T}
\newcommand{\de}{\delta_{\ee}}
\newcommand{\Phie}{\Phi_\ee}
\newcommand{\Thee}{\theta_\ee}

\renewcommand{\Re}{R_\ee}
\newcommand{\re}{r_\ee}

\newcommand{\cPrr}{\cP_{\re, \Re}}
\newcommand{\Cdn}{C_{\delta, \nu}}

\newtheorem{thm}{Theorem}
\newtheorem{lem}{Lemma}

\newtheorem{prop}{Proposition}
\newtheorem{defi}{Definition}

\newtheorem{rem}{Remark}

\title{Sharper dispersive estimates and asymptotics for a Boussinesq-type system with larger ill-prepared initial data}

\author{Fr\'ed\'eric Charve\footnote{Univ Paris Est Creteil, CNRS, LAMA, F-94010 Creteil, 2 Univ Gustave Eiffel, LAMA, F-77447 Marne-la-Vall\'ee, France. E-mail: frederic.charve@u-pec.fr}}

\date{}
\begin{document}

\maketitle

\begin{abstract} The aim of this article is to extend previous works about the asymptotics of an ill-prepared fast rotating, highly stratified incompressible Navier-Stokes system. Thanks to improved Strichartz and a priori estimates, we are able not only to cover a case which was unanswered in our previous work (allowing bigger ill-prepared initial data) but also to improve the convergence rates and reduce some assumptions on the initial data. In passing we also widen the range of some parameters and compare two methods to obtain dispersive estimates.
\end{abstract}

\section{Introduction}

\subsection{Geophysical fluids}

We refer to \cite{Chemin2, BMN5, FC, FCpochesLp, FCF1, FCPAA} for an introduction to the rotating and stratified Navier-Stokes system called the Primitive System (sometimes also called Primitive Equations) and to \cite{BeBo, Cushman, Sadourny, Pedlosky} for an in-depth presentation. This system aims at describing geophysical fluids located at the surface of the Earth (thus lying in a large physical domain) under the assumption that the vertical components of velocity or height are much smaller than their horizontal counterparts. Such a fluid is influenced by two concurrent "forces": first, the Coriolis force induced by the rotation of the Earth around its axis, and seconds the vertical stratification of the density induced by gravity.

This creates concurrent horizontal and vertical rigidities in the fluid motion and structure, and to measure their influence on the dynamics of the fluid, physicists introduced the Rossby and Froude numbers, namely $Ro$ and $Fr$. The smaller they are, the more influent are these two forces. Let us consider the Primitive System in the whole space, when both phenomena are of the same scale (that is we choose $Ro=\ee$ and $Fr=\ee F$ with $F>0$). In this article, to simplify, we call $\ee$ the Rossby number and $F$ the Froude number, and we write the system as follows:
\begin{equation}
\begin{cases}
\d_t \Ue +\ve\cdot \n \Ue -L \Ue +\frac{1}{\ee} \cA \Ue=\frac{1}{\ee} (-\n \Phie, 0),\\
\div \ve=0,\\
{\Ue}_{|t=0}=U_{0,\ee}.
\end{cases}
\label{PE}
\tag{$PE_\ee$}
\end{equation}
The unknowns are on one hand $\Ue =(\ve, \Thee)=(\ve^1, \ve^2, \ve^3, \Thee)$, where $\ve$ denotes the velocity of the fluid and $\Thee$ the scalar potential temperature (linked to the density, temperature and salinity), and on the other hand $\Phie$, which is called the geopotential and gathers the pressure term and the centrifugal force. The diffusion operator $L$ is defined by
$$
L\Ue \overset{\mbox{def}}{=} (\nu \D \ve, \nu' \D \Thee),
$$
where $\nu, \nu'>0$ denote the kinematic viscosity and thermal diffusivity (both will be considered as viscosities). The last term $\ee^{-1}\cA$ gathers the rotation and stratification effects and the matrix $\cA$ is defined by
$$
\cA \overset{\mbox{def}}{=}\left(
\begin{array}{llll}
0 & -1 & 0 & 0\\
1 & 0 & 0 & 0\\
0 & 0 & 0 & F^{-1}\\
0 & 0 & -F^{-1} & 0
\end{array}
\right).
$$

\subsection{Weak and strong solutions}

As emphasized in \cite{FCPAA}, the skew-symmetry of $\cA$ implies that any classical energy method like the ones used to study the Navier-Stokes system (based on $L^2$ or $H^s/\Hs $ inner products), will not "see" the penalized terms so we easily adapt the Leray and Fujita-Kato theorems. For any fixed $\ee>0$, there exist global-in-time weak solutions if $U_{0,\ee}\in L^2$, and a unique local-in-time strong solution when $U_{0,\ee} \in \dot{H}^\frac{1}{2}$ (global for an initial data whose norm is bounded by $c\nu_0$ for some small $c>0$). We also have at our disposal weak-strong uniqueness results and blow up criteria: denoting as $\Ue$ the unique strong solution of System \eqref{PE} provided by the Fujita-Kato theorem and defined on $[0,T]$ for all $0<T<T_\ee^*$, if the lifespan $T_\ee^*$ is finite then:
\begin{equation}
 \int_0^{T_\ee^*} \|\n \Ue(\tau)\|_{\dot{H}^\frac{1}{2}(\R^3)}^2 d\tau=\infty.
 \label{critereexpl}
\end{equation}
Moreover, if $\Uoe \in \Hs$ (for some $s\in]-\frac32, \frac32[$) then we can also propagate the regularity as done for the Navier-Stokes system. This is true wether $F=1$ (non-dispersive regime, we refer to \cite{Chemin2, FCF1}) or $F\neq 1$ (dispersive regime).

In the present work, for a fixed $F\neq 1$, our interest is to study the convergence (and obtain convergence rates) when $\ee$ goes to zero (that is for fast rotating and highly stratified systems) in the case of \emph{ill-posed large initial data} (see below for more details).

\subsection{The limit system, the QG/osc decomposition}
\label{sol}
The limit system as the small parameter $\ee$ goes to zero (which can first be formally identified then justified, see for example \cite{Chemin2, FC, FC2}) suggests a particular decomposition of the initial data and the solution, which helps studying the asymptotics.

Let us recall that this limit system is a transport-diffusion system coupled with a Biot-Savart inversion law and is called the 3D quasi-geostrophic system:
\begin{equation}
\begin{cases}
\d_t \tOm +\tvqg .\n \tOm -\G \tOm =0,\\
\tUqg =(\tvqg ,\tilde{\theta}_{QG})=(-\partial_2, \partial_1, 0, -F\partial_3) \DF^{-1} \tOm,
\end{cases}
\label{QG}\tag{$QG$}
\end{equation}
where we set $\DF=\d_1^2 +\d_2^2 +F^2 \d_3^2$, and the operator $\G$ is defined by:
$$
\G \overset{def}{=} \D \DF^{-1} (\nu \d_1^2 +\nu \d_2^2+ \nu' F^2 \d_3^2),
$$
The quantity $\tOm=\d_1 \tvqg ^2 -\d_2 \tvqg ^1 -F \d_3 \tilde{\theta}_{QG}$ is called the potential vorticity.
\begin{rem}
\sl{In general (that is when $F\neq 1$ and $\nu\neq \nu'$) $\G$ is a tricky non-local diffusion operator of order 2. We refer to \cite{FCestimLp, FCpochesLp} for an in-depth study of $\G$ (neither its Fourier kernel nor singular integral kernel have a constant sign and no classical result can be used).}
\end{rem}
Led by the limit system we first introduce for any 4-dimensional vectorfield $U=(v, \theta)$ its potential vorticity $\Om(U)$:
$$
\Om(U)\overset{def}{=} \d_1 v^2 -\d_2 v^1 -F\d_3 \theta,
$$
then its orthogonal decomposition into quasi-geostrophic and oscillating (or oscillatory) parts (quite similar to the Leray or Helmholtz decompositions, we refer to \cite{FC,FC5}):
\begin{equation}
U_{QG}=\cQ (U) \overset{def}{=} \left(
\begin{array}{c}
-\d_2\\
\d_1\\
0\\
-F\d_3
\end{array}
\right) \D_F^{-1} \Om (U), \quad \mbox{and} \quad U_{osc}=\cP (U) \overset{def}{=} U-U_{QG}.
\end{equation}
We will say that a vectorfield $U$ is quasi-geostrophic when $U=\cQ U$, and we refer to Proposition 2 in \cite{FCPAA} (and \cite{Chemin2, FC, FC2, FCpochesLp, FCF1}) for more properties of the associated orthogonal projectors $\cQ$ and $\cP$. In particular we can rewrite System \eqref{QG} into the following equivalent velocity formulation:
\begin{equation}
\begin{cases}
\d_t \tUqg +\cQ(\tvqg .\n \tUqg) -\G \tUqg =0,\\
\tUqg =\mathcal{Q} (\tUqg ),\\
{\tilde{U}_{QG|t=0}= \tUoqg,}
\end{cases}
\label{QG1}\tag{$QG$}
\end{equation}
Not only can we adapt the Leray and Fujita-Kato theorems to System \eqref{QG}, but this system also enjoys more "2D"-features as described in Theorem \ref{Th1} below: with additional regularity assumptions we obtain global existence without asking smallness of the initial data (the global existence and estimates from \cite{FC2} in the case of a $H^1$ initial data were improved in \cite{FC4, FCPAA}).
\begin{rem}
\sl{The notion of well-prepared/ill-prepared initial data is related to the closeness of the initial data to the limit quasi-geostrophic structure, in the sense that its oscillating part is small/going to zero/large/blowing up as $\ee$ goes to zero. In the present article we will be able to obtain global existence for more ill-prepared initial data (that is with larger oscillating parts) than in \cite{FCPAA}.}
\label{illprepared}
\end{rem}
Going back to System \eqref{PE}, the strategy is then to introduce $\Ome=\Om(\Ue)$, and study separately $\Uqg=\cQ(\Ue)$ and $\Uosc=\cP(\Ue)$. Let us give a quick view of the survey of the results from \cite{FC, FC2, FC3, FC4, FC5, FCestimLp, FCpochesLp, FCPAA} presented in Section 1.3 from \cite{FCPAA} (there coupled in each case with the specific use of Strichartz estimates). First, following the methods developped by Chemin, Desjardins, Gallagher and Grenier in \cite{CDGG, CDGG2, CDGGbook} (a general separation of the solution into the sum of one part close to the limit  and some oscillatory part which is small in some spaces thanks to Strichartz estimates obtained through non-stationnary phase arguments) we focussed not only on global existence but also on a study of the limit as $\ee$ goes to zero. In the series of works \cite{FC, FC2, FC3, FC4, FC5, FCPAA}, we obtained the asymptotics for \eqref{PE} when $\ee$ goes to zero (both in the general case $\nu\neq \nu'$ and in the particular case $\nu=\nu'$ where many tools are simplified) in various configurations of ill-prepared data and we reduced the initial regularity, considered large data or evanescent viscosities. Let us mention that in \cite{FC3, FCpochesLp}, we were interested in the vortex patches setting (inspired by \cite{Dutrifoy1, TH1}) and we had to study precisely the nonlocal diffusion operator $\G$ of the limit system (see \cite{FCestimLp}). Obtaining Strichartz estimates requires the study of the eigenvalues of the associated matrix (in the Fourier space) to the linearized system \eqref{systdisp}.
\\

In the same spirit, \cite{Dutrifoy2, VSN, VSNS1} are devoted to the rotating fluid system, \cite{Dutrifoy1} focusses on the non-viscous version of System \eqref{PE}, \cite{Scro3} studies the stratified system, and for the periodic case we refer to \cite{Scro, Scro2}. Let us also mention \cite{CMX} using similar dispersive methods for the 2D-quasi-geostrophic system. The methods in the non-dispersive case $F=1$ are rather different and there are only a few works (see \cite{Chemin2, Dragos4, FCF1}).
\\

Another series of articles \cite{GIMS, HS, IT2, IT5, KLT, KLT2, IMT, LT} is devoted to study the global existence of mild solutions only. In \cite{IT5} Iwabuchi and Takada study local solutions for the rotating fluid system. This work was extended in \cite{IT2} to global existence (with smallness assumptions for initial data but with large bound). In both articles they follow and extend the dispersive methods of Dutrifoy in \cite{Dutrifoy2}. In \cite{KLT} (and \cite{KLT2} for the non viscous case) Koh, Lee and Takada manage to widen the range for the parameters thanks to the use of the Littman theorem, which extends the classical stationnary phase method to the case where the hessian of the phase can degenerate, and allows to obtain slightly better dispersive and Strichartz estimates in the sense that the parameters are in a larger domain (we refer to Section \ref{Disp} for details). In \cite{IMT} Iwabuchi, Mahalov and Takada adapt this to System \eqref{PE} (only in the case $\nu=\nu'$) and in \cite{LT} Lee and Takada adapt the result in the case of stratification only (no rotationnal effects), also when $\nu=\nu'$.\\

More precisely, Iwabuchi, Mahalov and Takada obtained in \cite{IMT}  the following Strichartz estimates that we state with our notations:
\begin{prop} (\cite{IMT} Theorem 1.1 and Corollary 1.2)
\sl{Assume $F\neq 1$ and $\nu=\nu'$. If $r\in]2,4[$ and $p\in]2,\infty[ \cap [\frac1{2(\frac12-\frac1{r})}, \frac2{3(\frac12-\frac1{r})}]$, there exists a constant $C=C_{F,p,r}$ such that if $f$ solves the homogeneous system associated to \eqref{systdisp}, then
$$
\|f\|_{L^p(\R_+, L^r)} \leq C \frac{\ee^{\frac1{p} -\frac32(\frac12 -\frac1{r})}}{\nu^{\frac32(\frac12-\frac1{r})}} \|f_0\|_{L^2}.
$$
If $s\in]\frac12, \frac58]$, there exists a constant $C=C(F,s)$ such that:
$$
\|f\|_{L^4(\R_+, \dot{W}^{s,\frac6{1+2s}})} \leq C \frac{\ee^{\frac12(s-\frac12)}}{\nu^{\frac12(1-s)}} \|f_0\|_{\dot{H}^s}.
$$
}
\label{StriIMT}
\end{prop}
These estimates help them obtaining (when $\nu=\nu'$) through a fixed point argument the following global well-posedness results for initial data (independant of $\ee$) with small quasi-geostrophic part (see Theorems 1.3 and 1.5 in \cite{IMT}):
\begin{itemize}
 \item If $s\in]\frac12, \frac58]$, there exist $\delta_1,\delta_2>0$ (depending on $\nu,F,s$) such that for any $\ee>0$ and any initial data $U_0= U_{0,QG} + U_{0,osc}$ with $(U_{0,QG}, U_{0,osc})\in \dot{H}^\frac12 \times \dot{H}^s$ and
\begin{equation}
  \begin{cases}
   \|U_{0,QG}\|_{\dot{H}^\frac12} \leq \delta_1,\\
  \|U_{0,osc}\|_{\dot{H}^s} \leq \delta_2 \ee^{-\frac12(s-\frac12)},
 \end{cases}
 \label{IMTpetitesse}
\end{equation}
 there exists a unique global mild solution in $L^4(\R_+,\dot{W}^{\frac12, 3}(\R^3))$.
 \item There exists $\delta>0$ such that for any initial data $U_0= U_{0,QG} + U_{0,osc}\in \dot{H}^\frac12$ with $\|U_{0,QG}\|_{\dot{H}^\frac12} \leq \delta$, there exists $\ee_0>0$ such that for any $0<\ee<\ee_0$, System \eqref{PE} has a unique global mild solution in $\mathcal{C}(\R_+, \dot{H}^\frac12) \cap L^4(\R_+,\dot{W}^{\frac12, 3}(\R^3))$.
\end{itemize}

Motivated by this, we generalized in \cite{FCPAA} our results from \cite{FC2, FC3, FC4} and obtained full asymptotics (global existence, limit and convergence rates) for very large ill-prepared initial data (less regular, depending on $\ee$ and with a norm reaching the size of a negative power of $\ee$). We also generalized \cite{IMT} in the sense that we can consider initial data with \emph{large quasi-geostrophic part} (with low frequencies assumptions) and provide solutions in \emph{homogeneous} energy spaces $\dot{E}^s$ both in the particular case $\nu=\nu'$ \emph{and} in the general case $\nu\neq \nu'$. Our methods rely on the special structures and properties of the 3D quasi-geostrophic system (see Proposition 2 in \cite{FCPAA}) and the main ingredients were precise energy estimates using nonlocal derivation operators, coupled with improved Strichartz estimates, that we here write in the case $\nu=\nu'$:
\begin{prop} (\cite{FCPAA}, Proposition 47)
 \sl{Assume $F\neq 1$. For any $d\in \R$, $r>4$, $q\geq 1$ and
 $$
 \theta\in ]0,\frac{\frac12 -\frac1r}{1-\frac{4}{r}}[\cap]0,1], \quad p\in[1, \frac{4}{\theta (1-\frac{4}{r})}],
 $$ 
 there exists a constant $C_{F, p,\theta,r}$ such that if $f$ solves \eqref{systdisp} for initial data $f_0$ and external force $\Fe$ both with zero divergence and potential vorticity, then
 \begin{equation}
  \||D|^d f\|_{\tilde{L}_t^p\dot{B}_{r, q}^0} \leq C_{F,p,\theta,r}\frac{\ee^{\frac{\theta}{4}(1-\frac{4}{r})}}{\nu^{\frac{1}{p}-\frac{\theta}{4}(1-\frac{4}{r})}}  \times \left( \|f_0\|_{\dot{B}_{2, q}^\sigma} +\int_0^t  \|\Fe(\tau)\|_{\dot{B}_{2, q}^\sigma} d\tau \right),
 \end{equation}
 where $\sigma= d+\frac32-\frac{3}{r}-\frac{2}{p}+\frac{\theta}{2} (1-\frac{4}{r})$.
\label{Estimdispnuold}
 }
\end{prop}
In Theorems 6, 15 and 21 from \cite{FCPAA} we obtained asymptotics with initial oscillating parts of the size $\ee^{-\gamma}$ for some $\gamma<\frac{\delta}2$ (if $\frac12+\delta$ is the high frequency regularity of the initial part, $\delta<3/26$ if $\nu=\nu'$, $\delta \leq 1/2$ if $\nu\neq \nu'$). If only interested in global existence results, we expect to be able to consider initial oscillating parts of the size $\ee^{-\frac{\delta}2}$ as done in \cite{IMT}, but due to technical limitations we were not able to do this.

In the present article we will mainly focus on the case $F\neq 1$ and $\nu=\nu'$ (we will sometimes make remarks about the results in the case $\nu\neq \nu'$).

\subsection{Statement of the results}

We use the same notations as in \cite{FCPAA}, for $s\in \R$ and $T>0$ we define the space:
$$
\dot{E}_T^s=\mathcal{C}_T(\Hs (\R^3)) \cap L_T^2(\dot{H}^{s+1}(\R^3)),
$$
endowed with the following norm ($\nu=\nu'$):
$$
\|f\|_{\dot{E}_T^s}^2 \overset{def}{=}\|f\|_{L_T^\infty \Hs }^2+\nu \int_0^T \|f(\tau)\|_{\dot{H}^{s+1}}^2 d\tau.
$$
When $T=\infty$ we simply denote $\dot{E}^s$ and the corresponding norm is taken over $\R_+$ in time.
\\
Let us decompose our general ill-prepared initial data into its oscillating and quasi-geostrophic parts: $\Uoe =\Uoosc +\Uoqg$. While the $QG$-part will converge to some $\tUoqg$ (without any smallness condition), we allow the oscillating part to be very large (in terms of the Rossby number $\ee$). The aim of \cite{FCPAA} was to generalize \cite{FC2,FC3,FC4} and obtain global existence (and explicit convergence rates) for large initial data, with small extra regularity and the biggest possible blowing-up initial oscillatory part (as a negative power of $\ee$). Let us emphasize that one of the keys was to use as little as possible energy estimates for the oscillations $\We$ or $\Wet$ (which produce large terms such as $\mbox{exp}(\|\Uoosc\|^2)$ through the use of Gronwall estimates, and do not allow as large initial data as we wish). This forced us to obtain more flexible Strichartz estimates.
\\

Thanks to an improvement in our dispersive estimates allowing us to cover the range $r>2$ in the Strichartz estimates (instead of only $r>4$), we will be able to extend the work of \cite{FCPAA} and cover bigger initial data when only seeking global existence in $\dot{E}^\frac12$. Namely we obtain global existence for an initial oscillating part bounded by $c \ee^{-\frac{\delta}2}$ (where $c$ is a small constant) instead of $C \ee^{-\gamma}$ (with $\gamma<\frac{\delta}{2}$ and $C>0$ is free) which is the same as in \cite{KLT, IMT}. In addition our method provides the limit and convergence rates, which are not studied in \cite{KLT, IMT}. As a by-product, we manage to reduce the assumptions on the initial oscillating part, and widen the range for the parameter $\delta$ in Theorem 15 from \cite{FCPAA}. We also explicitely bound the $L^2L^\infty$-norm of $|D|^\beta(\Ue-\tUqg)$ (for some small $\beta$) by a power of $\ee$. 
\\

Before stating our result let us briefly recall the auxiliary systems involved in the statements and we refer to \cite{FCPAA} for more details (as usual we will systematically write, for $f:\R^3 \rightarrow \R^4$, $f\cdot \nabla f=\sum_{i=1}^3 f_i \d_i f$). Let us begin with rewriting the primitive system, projecting onto the divergence-free vectorfields ($\mathbb{P}$ denotes the Leray projector):
\begin{equation}
\begin{cases}
{\d_t \Ue -\nu \D \Ue+\frac{1}{\varepsilon} \mathbb{P} \mathcal{A} \Ue =-\mathbb{P}(\Ue.\nabla \Ue).} \\
{{\Ue}_{|t=0}=\Uoe.}
\end{cases}
\label{pe2}
\end{equation}
As in the previous works we also rewrite \eqref{QG} as follows (see \cite{FC2, FCPAA} for details):
\begin{equation}
\begin{cases}
{\d_t \tUqg -\nu \D \Ue+\frac{1}{\ee} \mathbb{P} \mathcal{A} \tUqg=-\mathbb{P}(\tUqg.\nabla \tUqg) +G^b,}\\
{\tilde{U}_{QG|t=0}= \tUoqg.}
\end{cases}
\label{QG3}\tag{$QG$}
\end{equation}
with $G^b$ the divergence-free and potential vorticity-free vectorfield defined as
\begin{equation}
G^b \overset{def}{=} \mathbb{P} \mathcal{P} (\tUqg. \nabla \tUqg).
\label{G}
\end{equation}
Next we consider the following linear system (whose aim is not only to provide oscillation but also to absorb the constant term $G^b$, we refer to Remark 11 in \cite{FCPAA} for details):
\begin{equation}
 \begin{cases}
  \d_t \We -\nu \D \We +\frac{1}{\varepsilon} \mathbb{P} \mathcal{A} \We = -G^b,\\
  {\We}_{|t=0}=\Uoosc
 \end{cases}
\label{We}
 \end{equation}
When $F\neq 1$ and $\nu=\nu'$, $\We$ is fully oscillatory (contrary to the general case $\nu\neq \nu'$ where we also have to deal with its QG-part) and satisfies System \eqref{systdisp} and therefore enjoys Strichartz-type estimates. We then define $\de=\Ue -\tUqg -\We$, which satisfies the following system:
\begin{equation}
 \begin{cases}
  \d_t \de -L\de +\frac{1}{\ee} \mathbb{P} \mathcal{A} \de = \Sum_{i=1}^8 F_i,\\
  {\de}_{|t=0}= \Uoqg-\tUoqg,
 \end{cases}
\label{GE}
\end{equation}
with:
\begin{equation}
\begin{cases}
 F_1 \overset{def}{=}-\mathbb{P}(\de \cdot \nabla \de), \quad F_2 \overset{def}{=}-\mathbb{P}(\de \cdot \nabla \tUqg), \quad F_3 \overset{def}{=}-\mathbb{P}(\tUqg \cdot \nabla \de),\\
 F_4 \overset{def}{=}-\mathbb{P}(\de \cdot \nabla \We),\quad F_5 \overset{def}{=}-\mathbb{P}(\We \cdot \nabla \de),\quad F_6 \overset{def}{=}-\mathbb{P}(\tUqg \cdot \nabla \We),\\
 F_7 \overset{def}{=}-\mathbb{P}(\We \cdot \nabla \tUqg),\quad F_8 \overset{def}{=}-\mathbb{P}(\We \cdot \nabla \We).
\end{cases}
 \label{f1f2}
\end{equation}
Let us state the following theorem, extending the global existence result for the limit system from \cite{FC2} thanks to the improved estimates from \cite{FC4} (here in the case $\nu=\nu'$):
\begin{thm} (\cite{FCPAA}, Theorem 14)
 \sl{Let $\delta\in]0,\frac12]$ and $\tUoqg\in H^{\frac12 + \delta}$ a quasi-geostrophic vectorfield (that is $\tUoqg=\cQ \tUoqg$). Then System \eqref{QG3} has a unique global solution in $E^{\frac12 + \delta}=\dot{E}^0\cap \dot{E}^{\frac12 + \delta}$ and there exists a constant $C=\Cdn>0$ such that for all $t\in \R_+$:
\begin{equation}
 \|\tUqg \|_{L_t^\infty H^{\frac12 + \delta}}^2+\nu \int_0^t \|\n \tUqg (\tau)\|_{H^{\frac12 + \delta}}^2 d\tau \leq \Cdn \|\tUoqg\|_{H^{\frac12 + \delta}}^2 \max(1, \|\tUoqg\|_{H^{\frac12 + \delta}}^{\frac{1}{\delta}}).
 \label{estimQG}
\end{equation}
\label{Th1}
 }
\end{thm}
We are now able to state our results. Let us begin with the following extension of Theorem 15 from \cite{FCPAA}:
\begin{thm}
 \sl{Assume $F\neq 1$ and $\nu=\nu'$. For any $\C_0\geq 1$, $\delta\in]0,\frac14[$, $\gamma\in]0,\frac{\delta}{2}[$ and any $\aa_0>0$, if we define $\eta_0>0$ such that
 \begin{equation}
  \gamma=(1-2\eta_0)\frac{\delta}{2} \quad (\mbox{that is } \eta_0=\frac12(1-\frac{2\gamma}{\delta})),
 \label{defdeltaeta}
 \end{equation}
 there exist $\ee_0,\B_0>0$ (all of them depending on $F, \nu,\Co, \delta, \gamma, \alpha_0$) such that for all $\ee\in]0,\ee_0]$ and all divergence-free initial data $\Uoe= \Uoqg + \Uoosc$ satisfying:
 \begin{enumerate}
  \item There exists a quasi-geostrophic vectorfield $\tUoqg\in H^{\frac12 + \delta}$ such that
  \begin{equation}
   \begin{cases}
    \|\Uoqg-\tUoqg\|_{H^{\frac12 + \delta}}\leq \C_0 \ee^{\aa_0},\\
    \|\tUoqg\|_{H^{\frac12 + \delta}}\leq \C_0.
   \end{cases}
  \end{equation}
\item $\Uoosc\in \dot{H}^{\frac12} \cap \dot{H}^{\frac12 + \delta}$ with $\|\Uoosc\|_{\dot{H}^{\frac12 + \delta}} \leq \C_0 \ee^{-\gamma}$,
 \end{enumerate}
then the following results are true:
\begin{enumerate}
 \item System \eqref{PE} has a unique global solution $\Ue \in \dot{E^s}$ for all $s\in[\frac12, \frac12+2\eta_0 \delta[$, and if we define
\begin{itemize}
 \item $\tUqg$ as the unique global solution of \eqref{QG3} in $\dot{E}^0 \cap \dot{E}^{\frac12 +\delta}$,
\item $\We$ as the unique global solution of \eqref{We} in $\dot{E}^{\frac12} \cap \dot{E}^{\frac12 +\delta}$,
\item $\de=\Ue-\tUqg-\We$,
\end{itemize}
then for all $s\in[\frac12, \frac12+2\eta_0 \delta[$
$$
\|\de\|_{\dot{E}^s}\leq \B_0 \ee^{\min\big(\aa_0, \frac{\delta}2 -\gamma, \frac{\delta}2-\gamma +\frac12 (\frac12-s)\big)}\leq \Bo \ee^{\min\big(\aa_0, \eta_0 \delta-\frac12 (s-\frac12) \big)}.
$$
\item Moreover if, in addition, we ask $\|\Uoosc\|_{\dot{H}^{\frac12+ c\delta} \cap \dot{H}^{\frac12 + \delta}} \leq \C_0 \ee^{-\gamma}$ (for some fixed $c \in ]0,1[$ close to 1) then we can get rid of the oscillations: for all $\eta\in]0,2\eta_0[$, for all $\eta'\in]0,\min(\eta, c)[$,
$$
\big\||D|^{\eta'\delta}(\Ue-\tUqg)\big\|_{L^2 L^\infty} \leq \Bo \ee^{\min\big(\aa_0, (\eta_0-\frac{\eta}2)\delta\big)}.
$$
\item Finally, if we ask more low-frequency regularity on the initial oscillating part, that is $\Uoosc\in \dot{H}^{\frac12-\delta} \cap \dot{H}^{\frac12 + \delta}$ with the same condition $\|\Uoosc\|_{\dot{H}^{\frac12+ c\delta} \cap \dot{H}^{\frac12 + \delta}} \leq \C_0 \ee^{-\gamma}$ (for some fixed $c \in ]0,1[$ close to 1), then the estimates
$$
\|\de\|_{\dot{E}^s}\leq \B_0 \ee^{\min\big(\aa_0, \frac{\delta}2 -\gamma, \frac{\delta}2-\gamma +\frac12 (\frac12-s)\big)}\leq \B_0 \ee^{\min\left(\aa_0, (\eta_0-\frac{\eta}2)\delta\right)},
$$
remain valid for all $s\in[\frac12 -\eta \delta, \frac12+\eta \delta]$ (with $0<\eta <\min\big(1-2\frac{\gamma}{\delta}, \frac1{3\delta}-1, \frac1{2\delta}-2\big)$) and we also obtain that:
$$
\|\Ue-\tUqg\|_{L^2 L^\infty} \leq \B_0 \ee^{\min(\aa_0, k\eta_0\delta)},
$$
for any $k<1$ as close to $1$ as we wish.
\end{enumerate}
\label{Th2}
 }
\end{thm}
\begin{rem}
 \sl{\begin{enumerate}
      \item The first improvement is that we only ask $\|\Uoosc\|_{\dot{H}^{\frac12 +\delta}} \leq \C_0 \ee^{-\gamma}$ and in points 2 and 3, we only ask $\|\Uoosc\|_{\dot{H}^{\frac12 +c \delta} \cap \dot{H}^{\frac12 + \delta}} \leq \C_0 \ee^{-\gamma}$ for some $c<1$ close to 1 (instead of bounding respectively the $\dot{H}^\frac12 \cap \dot{H}^{\frac12 + \delta}$ and $\dot{H}^{\frac12-\delta} \cap \dot{H}^{\frac12 + \delta}$ norms in \cite{FCPAA}).
      \item The $L^2 L^\infty$ estimate from point 3 is slightly more precise and better than in \cite{FCPAA} (where the convergence rate was $\ee^{\min(\aa_0, \eta_0\frac{\delta}2)}$).
      \item A particular case of the estimates from point 2 says that (taking $\eta'=\eta/2$) for all $\eta\in]0,2\min(\eta_0,c)[$,
$$
\big\||D|^{\eta \frac{\delta}2}(\Ue-\tUqg)\big\|_{L^2 L^\infty} \leq \Bo \ee^{\min\big(\aa_0, (\eta_0-\frac{\eta}2)\delta\big)}.
$$
     \item Our new Strichartz estimates allow us to push a little further the upper bound for $\delta>0$. More precisely, $\delta< \frac14$ (In Theorem 15 from \cite{FCPAA} the result is given for $0<\delta\leq \frac1{10}$ for simplicity but the real bound was $0<\delta< \frac3{26}$, and $\delta\leq\frac1{8}$ in \cite{IMT}). Had we not used the Littman theorem we would end up with the same result but for $\delta<\frac18$ (see Section \ref{RemLitt}).
      \item In \cite{FCPAA} we managed to avoid using energy estimates for the oscillations $\We$ except on two terms. In the present article, thanks to our extended Strichartz estimates we now do not use them at all. 
      \item In the general case $\nu\neq\nu'$, the same improvement can be performed on Proposition 51 from \cite{FCPAA} (with $r\geq 2$, $p \leq 4/(1-\frac2{r})$ and slightly better coefficients) but there will be no other improvement as the convergence rate is limited by terms coming from the various truncations. Similar estimates on $\big\||D|^{\eta \frac{\delta}2}(\Ue-\tUqg)\big\|_{L^2 L^\infty}$ as in point 2 can be proved without difficulty.
      \end{enumerate}
}
\end{rem}
Let us now state the second result of this article, which answers Remarks 19 and 32 from \cite{FCPAA} (and which is the first motivation in the improvement of the Strichartz estimates):
\begin{thm}
 \sl{Assume $F\neq 1$ and $\nu=\nu'$. For any $\C_0\geq 1$, $\delta\in]0,\frac14[$, and any $\aa_0>0$, there exist $\ee_0, m_0>0$ (both of them small and depending on $F, \nu,\Co, \delta, \alpha_0$) such that for all $\ee\in]0,\ee_0]$ and all divergence-free initial data $\Uoe= \Uoqg + \Uoosc$ with with $\Uoqg$ (and $\tUoqg$) as in Theorem \ref{Th2} and $\Uoosc\in \dot{H}^{\frac12} \cap \dot{H}^{\frac12 + \delta}$ with
 $$
 \|\Uoosc\|_{\dot{H}^{\frac12 + \delta}} \leq m_0 \ee^{-\frac{\delta}2},
 $$
then System \eqref{PE} has a unique global solution $\Ue \in \dot{E}^\frac12$.
\\
Moreover if we assume there exists a function $m(\ee)$ such that $m(\ee)\underset{\ee\rightarrow 0}{\longrightarrow} 0$ and:
\begin{equation}
\|\Uoosc\|_{\dot{H}^{\frac12 + \delta}} \leq m(\ee) \ee^{-\frac{\delta}2},
\end{equation}
then in addition, if we define $\tUqg$, $\We$ and $\de$ as previoulsy, there exists $\B_0$ (depending on $F, \nu,\Co, \delta, \alpha_0$) such that we have
\begin{equation}
 \|\de\|_{\dot{E}^\frac12}\leq \B_0 \max(\ee^{\aa_0}, m(\ee)).
 \label{estimThdelta}
\end{equation}
\label{Th3}
 }
\end{thm}

\begin{rem}
 \sl{\begin{enumerate}
      \item Actually, in the proof of the first case, we obtain that for some $C\geq 1$, $\|\de\|_{\dot{E}^\frac12}\leq \B_0 \min(\ee^{\aa_0},m_0) \leq \frac{\nu}{C}$.
      \item Compared to \cite{KLT, IMT} we now reach the same size of $\ee^{-\frac{\delta}{2}}$ when we only inquire the global well-posedness in $\dot{E}^\frac12$. As outlined in Remark 17 in \cite{FCPAA}, we do not have any smallness condition on the $QG$-part (contrary to \cite{KLT, IMT}) and we also provide a convergence rate. This totally answers Remarks 19 and 32 from \cite{FCPAA} as outlined in the following remark.
      \item As we consider energy solutions, we still require an extra low frequency regularity for the initial oscillating part, and extra regularity for the initial $QG$-part compared to \cite{IMT} (see remark 18 in \cite{FCPAA}).
      \end{enumerate}
     }
\end{rem}

\begin{rem}
 \sl{We can mirror Remark 32 from \cite{FCPAA} as we prove that:
 \begin{itemize}
  \item if $\Uoosc\in \dot{H}^\frac12 \cap \dot{H}^{\frac12+\delta}$ and $\|\Uoosc\|_{\dot{H}^{\frac12+\delta}} \ee^{\frac{\delta}2}\leq m_0$, with $m_0$ small enough, we obtain global well-posedness in $\dot{E}^\frac12$ for $\ee\leq \ee_0$,
  \item if $\|\Uoosc\|_{\dot{H}^{\frac12+\delta}} \ee^{\frac{\delta}2}\underset{\ee\rightarrow 0}{\rightarrow} 0$, then in addition $\|\de\|_{\dot{E}^\frac12}$ goes to zero (even if $\|\Ue\|_{\dot{E}^\frac12}$ may blow up as $\ee$ goes to zero),
  \item if $\|\Uoosc\|_{\dot{H}^{\frac12+\delta}} \leq \Co \ee^{-\gamma}$, with $\gamma<\frac{\delta}{2}$, then we get a convergence rate of $\|\de\|_{\dot{E}^s}$ (for all $s\in [\frac12, \frac12+2\eta_0 \delta[$) as a positive power of $\ee$.
  \item if, in addition to the previous case, we ask $\|\Uoosc\|_{\dot{H}^{\frac12+c\delta}\cap \dot{H}^{\frac12+\delta}} \leq \Co \ee^{-\gamma}$ (for some $c\in]0,1[$), then we can bound $\||D|^\beta(\Ue-\tUqg)\|_{L^2 L^\infty}$ (for some small $\beta$) by a positive power of $\ee$.
  \item if $\Uoosc\in \dot{H}^{\frac12-\delta} \cap \dot{H}^{\frac12+\delta}$ and $\|\Uoosc\|_{\dot{H}^{\frac12+c\delta}\cap \dot{H}^{\frac12+\delta}} \leq \Co \ee^{-\gamma}$ (for $c\in]0,1[$), with $\gamma<\frac{\delta}{2}$, then in addition $\|\Ue-\tUqg\|_{L^2 L^\infty}$ is bounded by a positive power of $\ee$.
 \end{itemize}
}
 \label{RqIMT}
\end{rem}
This article is structured as follows: we begin with the new Strichartz estimates, then we improve the a priori estimates and finally we prove Theorems \ref{Th2} and \ref{Th3}. For the sake of conciseness we will only focus on what is new and often refer to \cite{FCPAA} for precisions. For the same reason we will give minimal notational informations about the Littlewood-Paley decomposition and refer to \cite{Dbook} for an in-depth description. We end the article with a remark about the fact that these new results can be also obtained without resorting to the Littman theorem (we would then have a slightly smaller range for the parameter $\delta$).

\section{Dispersion and Strichartz estimates}

\label{Disp}

Consider the following system where both $\Fe$ and $f_0$ are divergence-free (in the case $\nu=\nu'$, we have $L=\nu\Delta$):
\begin{equation}
\begin{cases}
\d_t f-(\nu\Delta-\frac{1}{\ee} \mathbb{P} \cA) f=\Fe,\\
f_{|t=0}=f_0.
\end{cases}
\label{systdisp}
\end{equation}
If we apply the Fourier transform, the equation becomes (see \cite{FC, FCPAA} for precisions):
$$
\d_t \hat{f}- \mathbb{B}(\xi, \ee)\hat{f}=\hat{\Fe},
$$
where
$$\mathbb{B}(\xi, \ee)= \hat{L-\frac{1}{\ee} \mathbb{P} \cA} =\left(
\begin{array}{cccc}
\displaystyle{-\nu|\xi|^2+\frac{\xi_1\xi_2}{\ee
  |\xi|^2}} & \displaystyle{\frac{\xi_2^2+\xi_3^2}{\ee
  |\xi|^2}} & 0 & \displaystyle{\frac{\xi_1\xi_3}{\ee F |\xi|^2}}\\
\displaystyle{-\frac{\xi_1^2+\xi_3^2}{\ee
  |\xi|^2}} & \displaystyle{-\nu|\xi|^2-\frac{\xi_1\xi_2}{\ee
  |\xi|^2}} & 0 & \displaystyle{\frac{\xi_2\xi_3}{\ee F |\xi|^2}}\\
\displaystyle{\frac{\xi_2\xi_3}{\ee |\xi|^2}} &
\displaystyle{-\frac{\xi_1\xi_3}{\ee
  |\xi|^2}} & \displaystyle{-\nu |\xi|^2} & \displaystyle{-\frac{\xi_1^2+\xi_2^2}{\ee F
  |\xi|^2}}\\
0 & 0 & \displaystyle{\frac{1}{\ee F}} &
\displaystyle{-\nu|\xi|^2}
\end{array}
\right).
$$

\subsection{Eigenvalues, projectors}

We begin with the eigenvalues and eigenvectors of matrix $\mathbb{B}(\xi, \ee)$. When $\nu=\nu'$, one of the four eigenvalues is discarded (as its associated eigenvector is not "divergence-free"), the second one will correspond to the QG-part, and the last two will correspond to the oscillating part.

We recall that when $\nu=\nu'$, the operator $\G$ reduces to $\nu \Delta$ and we refer to \cite{FCPAA} (Proposition 43) for details in the general case (no assumption for $\nu,\nu'$) about the following proposition.

\begin{prop}
\label{estimvp}
\sl{
If $\nu=\nu'$ the matrix $\mathbb{B}(\xi, \ee) = \widehat{L-\frac{1}{\ee} \mathbb{P} \mathcal{A}}$ is diagonalizable when restricted to the subspace $(\xi_1, \xi_2, \xi_3,0)^\bot$ and its eigenvalues are (where we denote $|\xi|_F^2 = \xi_1^2 + \xi_2^2 + F^2\xi_3^3$):
\begin{equation}
\label{vp}
\begin{cases}
\vspace{0.2cm} \mu = -\nu |\xi|^2,\\
\vspace{0.2cm} \lambda = -\nu |\xi|^2 +i\frac{|\xi|_F}{\ee F|\xi|},\\
\vspace{0.2cm} \overline{\lambda} = -\nu |\xi|^2 -i\frac{|\xi|_F}{\ee F|\xi|},
\end{cases}
\end{equation}
The projectors $\mathcal{P}_i(\xi, \ee)$ (we keep the indices $i\in\{2,3,4\}$ as in the general case) onto the eigenspaces corresponding to $\mu$, $\lambda$ and $\overline{\lambda}$, are of norm $1$ and mutually orthogonal projectors.

Moreover, setting $\mathbb{P}_i(u)= \mathcal{F}^{-1}\big(\mathcal{P}_i (\xi, \ee)(\widehat{u}(\xi))\big)$, in the case $\nu=\nu'$ they exactly fit the $QG/osc$ decomposition (for divergence-free vectorfields) in the sense that:
$$
\cP=\mathbb{P}_{3+4}=\mathbb{P}_3+\mathbb{P}_4 \quad \mbox{and }\cQ=\mathbb{P}_2.
$$
\label{remnunu}
 }
\end{prop}

\subsection{Strichartz estimates in the case $\nu=\nu'$}

The main result of this section is stated as follows (we recall the definition of Chemin-Lerner-Besov spaces as well as some properties in the appendix):

\begin{prop}
 \sl{There exists a constant $C_F>0$, such that for any $d\in \R$, $r\geq 2$, $q\geq 1$, $\theta\in[0,1]$ and $p\in[1, \frac{4}{\theta (1-\frac{2}{r})}]$, if $f$ solves \eqref{systdisp} for initial data $f_0$ and external force $\Fe$ both with zero divergence and potential vorticity, then
 \begin{equation}
  \||D|^d f\|_{\tilde{L}_t^p\dot{B}_{r, q}^0} \leq \frac{C_{F, p,\theta,r}}{\nu^{\frac{1}{p}-\frac{\theta}{4}(1-\frac{2}{r})}} \ee^{\frac{\theta}{4}(1-\frac{2}{r})} \times \left( \|f_0\|_{\dot{B}_{2, q}^\sigma} +\int_0^t  \|\Fe(\tau)\|_{\dot{B}_{2, q}^\sigma} d\tau \right),
 \end{equation}
 where
 $$
 \begin{cases}
  \sigma= d+\frac32-\frac{3}{r}-\frac{2}{p}+\frac{\theta}{2} (1-\frac{2}{r}),\\
  C_{F, p,\theta,r}=(C_F)^{\frac12-\frac1{r}}\left[4\Big(\frac{1}{p}-\frac{\theta}{4}(1-\frac{2}{r})\Big)\right]^{\frac{1}{p}-\frac{\theta}{4}(1-\frac{2}{r})}.
 \end{cases}
 $$
\label{Estimdispnu}
 }
\end{prop}

\begin{rem}
\sl{The estimates from \cite{IMT} (Theorem 1.1 and Corollary 1.2) recalled in Proposition \ref{StriIMT} are particular cases of the previous estimates:
\begin{itemize}
 \item To get the first one we simply need to take $\theta=\frac2{1-\frac2{r}}(\frac3{r}+\frac2{p}-\frac32)$ which is possible if and only if $\frac2{p} \in [3(\frac12 -\frac1{r}), 4(\frac12 -\frac1{r})]$.
 \item For the second one, we can choose $\theta=\frac3{1-s}(s-\frac12)$, which is possible if and only if $s\in [\frac12, \frac58]$.
\end{itemize}
}
\end{rem}

\textbf{Proof of Proposition \ref{Estimdispnu}: } Let us assume that $\Fe=0$ (if not, just repeat the following argument to the Duhamel term). As explained in \cite{FCPAA}, as $f_0$ is divergence-free and with zero potential vorticity we have:
$$
f_0=\mathbb{P} f_0=\cP \mathbb{P} f_0= \mathbb{P}_{3+4} \mathbb{P} f_0 =\mathbb{P}_{3+4} f_0,
$$
and thanks to the orthogonal decomposition
$$
f(t)=\mathcal{F}^{-1}\left(e^{-\nu t |\xi|^2 +i\frac{t}{\ee} \frac{|\xi|_F}{F|\xi|}} \mathcal{P}_3(\xi,\ee) \hat{f_0}(\xi) +e^{-\nu t |\xi|^2 -i\frac{t}{\ee} \frac{|\xi|_F}{F|\xi|}} \mathcal{P}_4(\xi,\ee) \hat{f_0}(\xi) \right),
$$
so that we write, in order to simplify:
$$
f(t)=\mathcal{F}^{-1}\left(e^{-\nu t |\xi|^2 +i\frac{t}{\ee} \frac{|\xi|_F}{F|\xi|}} \hat{f_0}(\xi) \right),
$$
If $\varphi$ is the truncation function involved in the Littlewood-Paley decomposition, we denote by $\varphi_1$ another smooth truncation function, with support in a slightly larger annulus than $\mbox{supp }\varphi$ (say for example the annulus centered at zero and of radii $\frac12$ and $3$) and equal to $1$ on $\mbox{supp }\varphi$. Let $\cB$ be the set:
$$
\cB\overset{def}{=}\{\psi \in \cC_0^\infty (\R_+\times \R^3, \R), \quad \|\psi\|_{L^{\bar{p}}(\R_+, L^{\bar{r}}(\R^3))}\leq 1\},
$$
then we follow the same classical steps, except that (as done in \cite{Dutrifoy2}) we will avoid keeping the $\psi$ functions into some convolution product which was the way we did in $(A.88)$ from \cite{FCPAA} (this forces the other term to be estimated in $L^{\frac{r}2}$ and this is where the limitation $r\geq4$ comes from): for any $j\in \Z$ and $r \geq 1$:
\begin{multline}
 \|\ddj f\|_{L^p L^r}= \sup_{\psi \in \cB} \int_0^\infty \int_{\R^3} \ddj f(t,x) \psi(t,x) dx dt\\
 =C \sup_{\psi \in \cB} \int_0^\infty \int_{\R^3} e^{-\nu t|\xi|^2+i\frac{t}{\ee}\frac{|\xi|_F}{F|\xi|}} \widehat{\ddj f_0}(\xi) \varphi_1(2^{-j} \xi) \hat{\psi}(t,\xi) d\xi dt\\
 \leq C \|\ddj f_0\|_{L^2} \sup_{\psi \in \cB} \left[\int_{\R^3} \int_0^\infty \int_0^\infty e^{-\nu (t+t')|\xi|^2+i\frac{t-t'}{\ee}\frac{|\xi|_F}{F|\xi|}} \varphi_1(2^{-j} \xi)^2 \hat{\psi}(t,\xi) \overline{\hat{\psi}(t',\xi)} dt dt'd\xi\right]^{\frac12}\\
 \leq C \sup_{\psi \in \cB} \|\ddj f_0\|_{L^2} \left[\int_0^\infty \int_0^\infty \|L_j(\frac{t-t'}{\ee})\psi(t,.)\|_{L^r} \|e^{\nu(t+t')\Delta} \varphi_1(2^{-j}D) \overline{\psi(t',.)}\|_{L^{\bar{r}}} dtdt'\right]^{\frac12},
\label{estimTT1}
 \end{multline}
with $L_j(\sigma)$ defined as follows:
\begin{equation}
 L_j(\sigma)g=\int_{\R^3} e^{ix\cdot \xi +i \sigma \frac{|\xi|_F}{F|\xi|}} \varphi_1(2^{-j}|\xi|) \hat{g}(\xi) d\xi =K_j(\sigma)*g,
\label{defLj}
\end{equation}
where
\begin{equation}
K_j(\sigma)(x)=\int_{\R^3} e^{ix\cdot \xi +i \sigma \frac{|\xi|_F}{F|\xi|}} \varphi_1(2^{-j}|\xi|) d\xi.
 \label{defKj}
\end{equation}
Thanks to the frequency truncation (remember that $\mbox{supp }\varphi_1 \subset \cC(0,\frac12,3)$) and the classical estimates for the heat kernel in this case (we refer for example to Lemma 2.4 from \cite{Dbook}) we easily get that:
\begin{equation}
 \|e^{\nu(t+t')\Delta} \varphi_1(2^{-j}D) \overline{\psi(t',.)}\|_{L^{\bar{r}}} \leq C'e^{-\frac{\nu}4 (t+t')2^{2j}} \|\psi(t')\|_{L^{\bar{r}}}.
\label{Estimphi}
 \end{equation}
For the other term we will take advantage of the Riesz-Thorin and Littman theorems. The latter extends the classical stationnary phase results known for a nondegenerate Hessian (we refer to Proposition 5 and Theorem 1 from \cite{Stein2} chapter VIII sections 2 and 3) into the case of a Hessian with $k$ nonzero eigenvalues (we refer to \cite{Litt} and 5.8 in \cite{Stein} chapter VIII section B). For the convenience of the reader we recall this result with our notations:
\begin{thm} (Littman \cite{Litt, Stein2})
 \sl{Assume that $\psi:\R^n \rightarrow \R$ is a smooth function compactly supported in $K$ and $\phi:\R^n \rightarrow \R$ is a smooth function such that for any $\xi \in K$, the Hessian $D^2 \phi (\xi)$ has at least $k$ nonzero eigenvalues. Then there exists a constant $A$ such that for any $\lambda\in \R$ and $x\in\R^n$,
 $$
 |\int_{\R^n} e^{ix\cdot \xi +i \lambda \phi(\xi)} \psi(\xi) d\xi| \leq A\sqrt{|x|^2+\lambda^2}^{-\frac{k}2} \leq A |\lambda|^{-\frac{k}2}.
 $$
 }
 \label{ThLitt}
\end{thm}

\begin{rem}
 \sl{As $\psi$ is compactly supported, using that:
 $$
 |\int_{\R^n} e^{ix\cdot \xi +i \lambda \phi(\xi)} \psi(\xi) d\xi| \leq \|\psi\|_{L^1},
 $$
 there exists a constant $C$ such that:
 \begin{equation}
   |\int_{\R^n} e^{ix\cdot \xi +i \lambda \phi(\xi)} \psi(\xi) d\xi| \leq C \min(1,|\lambda|^{-\frac{k}2}),
\label{ThLitt2}
 \end{equation}
 which is the way we will use the Littman theorem in what follows.}
\end{rem}

\begin{rem}
 \sl{Note that thanks to the previous result we do not need anymore to use a vertical Littlewood-Paley decomposition in order to get rid of the singularity and use simple non-stationnary argument (as done in \cite{CDGG, CDGGbook} and in many other works, for example \cite{FC2, FC3, FCPAA}). This vertical decomposition required a summation in $k$ that was limiting the range of $r,\theta$ as explained in Remark \ref{plusdek}.}
\end{rem}

Even if we do not have to deal with $\psi(t)*\psi(t')$ anymore, if we classically write the Young estimate:
$$
\|L_j(\sigma)g\|_{L^r}\leq \|K_j(\sigma)\|_{L^s}\|g\|_{L^{\bar{r}}},
$$
we can only choose $s=\frac{r}2$ and as the kernel $K_j(\sigma)$ is estimated in $L^q$ for $q\geq 2$ we would still need $r\geq 4$. This problem is avoided dealing with $L_j(\sigma) g$ through the Riesz-Thorin theorem (as done in \cite{Dutrifoy2, KLT, IMT}). First, considering $L_j$ instead of $K_j$, the $L^2$ estimates will not feature any derivative (no $2^j$) as thanks to the Plancherel formula there exists a constant $C$ (only depending on $\varphi_1$) such that :
\begin{equation}
 \|L_j(\sigma) g\|_{L^2} \leq C \|g\|_{L^2}.
 \label{EstimLjL2}
\end{equation}
Next, using the Young estimate:
\begin{equation}
 \|L_j(\sigma)g\|_{L^\infty}\leq \|K_j(\sigma)\|_{L^\infty}\|g\|_{L^1},
 \label{EstimLjL1}
\end{equation}
and thanks to the change of variable $\xi=2^j \eta$, we get $K_j(\sigma)(x)=2^{3j} K_0(\sigma)(2^j x)$ and:
\begin{equation}
 \|K_j(\sigma)\|_{L^\infty}\leq 2^{3j}\|K_0(\sigma)\|_{L^\infty},
\label{Chgtvar}
 \end{equation}
so that we are reduced to study with the help of the Littman theorem (as is done in \cite{KLT, IMT, LT}) the $L^\infty$-norm of $K_0$ (see \eqref{defKj} for the definition) that we will write as follows when denoting $b_F(\xi)=\frac{|\xi|_F}{F|\xi|}$:
\begin{equation}
 K_0(\sigma)(x)=\int_{\R^3} e^{ix\cdot \xi +i \sigma b_F(\xi)} \varphi_1(|\xi|) d\xi,
\end{equation}
Let us give some details about the strategy that consists in frequency truncations in order to decompose the frequency space into zones where we precisely know how many nonzero eigenvalues are featured by the Hessian $D^2 b_F$ which writes:
\begin{equation}
D^2 b_F(\xi)=\frac{1-F^2}{|\xi|^3|\xi|_F} \left(
 \begin{array}{ccc}
  \xi_3^2 (1-\xi_1^2 A_F(\xi)) & -\xi_1 \xi_2 \xi_3^2 A_F(\xi) & \xi_1 \xi_3 (2-\xi_3^2 B_F(\xi))\\
  -\xi_1 \xi_2 \xi_3^2 A_F(\xi) & \xi_3^2 (1-\xi_2^2 A_F(\xi)) & \xi_2 \xi_3 (2-\xi_3^2 B_F(\xi))\\
  \xi_1 \xi_3 (2-\xi_3^2 B_F(\xi)) & \xi_2 \xi_3 (2-\xi_3^2 B_F(\xi)) & -(\xi_1^2+\xi_2^2) (1-\xi_3^2 B_F(\xi))
 \end{array} \right),
\end{equation}
where we denote
$$
A_F(\xi)=\frac3{|\xi|^2}+\frac1{|\xi|_F^2} \quad \mbox{and } B_F(\xi)=\frac3{|\xi|^2}+\frac{F^2}{|\xi|_F^2}.
$$
Instead of reasoning with the determinant (as done in \cite{KLT, IMT}) let us simply clearly compute the eigenvalues of this matrix:
\begin{multline}
 \frac{1-F^2}{|\xi|^3|\xi|_F} \left\{\xi_3^2, \frac12\left(\xi_3^2-(\xi_1^2 +\xi_2^2)\frac{|\xi|^2}{|\xi|_F^2} \pm \sqrt{(\xi_3^2-(\xi_1^2 +\xi_2^2)\frac{|\xi|^2}{|\xi|_F^2})^2+4\xi_3^2(\xi_1^2+\xi_2^2)}\right) \right\}\\
 =\frac{1-F^2}{|\xi|^3|\xi|_F} \left\{\xi_3^2, \frac12\left(\frac{F^2\xi_3^4-|\xi_h|^4}{|\xi|_F^2} \pm \sqrt{(\frac{F^2\xi_3^4-|\xi_h|^4}{|\xi|_F^2})^2+4\xi_3^2(\xi_1^2+\xi_2^2)}\right) \right\}.
\end{multline}
One of the last two eigenvalues is equal to zero if and only if $\xi_3=0$ or $\xi_h=0$. More precisely, if $\xi_3=0$ and $\xi_h\neq 0$ then only one out of the three eigenvalues is nonzero, if $\xi_3\neq 0$ and $\xi_h=0$ then two of them are nonzero. When $\xi_3\neq 0$ and $\xi_h\neq 0$ then the three of them are nonzero.
\begin{rem}
 \sl{It is interesting to compare this with the case of the rotating fluids system where the eigenvalues of the hessian of the phase $b(\xi)=\frac{\xi_3}{|\xi|}$ are given by
$$
\frac1{|\xi|^3} \left\{-\xi_3, \xi_3\pm\sqrt{4|\xi_h|^2+\xi_3^2} \right\},
$$
and when $\xi_3=0$ and $\xi_h\neq 0$ we get two nonzero eigenvalues instead of only one in our case. This explains the better improvement (in terms of the power of $\ee$) given by the use of the Littman theorem in the case of the rotating fluids (see \cite{KLT}) compared to the case of the primitive system (see \cite{IMT}).
 }
 \label{remVPRF}
\end{rem}
Once this precision has been given, as in \cite{KLT, IMT} we will simply decompose as follows:
$$
K_0(\sigma)(x)=\int_{\R^3} e^{ix\cdot \xi +i \frac{\sigma}{F} b_F(\xi)} \Big(\chi(\frac{|\xi_3|}{k_0}) +\big(1-\chi(\frac{|\xi_3|}{k_0})\big) \Big) \Big(\chi(\frac{|\xi_h|}{k_0}) +\big(1-\chi(\frac{|\xi_h|}{k_0})\big) \Big) \varphi_1(|\xi|) d\xi,
$$
and when $k_0>0$ is chosen small enough (As $\mbox{supp }\varphi_1=\cC (0,\frac12, 3)$ and $\mbox{supp }\chi =B(0,\frac43)$, it suffices to ask $k_0<\frac3{8\sqrt{2}}$), we can get rid of the term involving small horizontal and vertical frequencies and write:
$$
K_0(\sigma)= \sum_{j=1}^3 K_{0,j}(\sigma),
$$
with
$$
\begin{cases}
\vspace{0.1cm}
 \displaystyle{K_{0,1}(\sigma)(x)= \int_{\R^3} e^{ix\cdot \xi +i \frac{\sigma}{F} b_F(\xi)} \chi(\frac{|\xi_3|}{k_0}) \big(1-\chi(\frac{|\xi_h|}{k_0})\big) \varphi_1(|\xi|) d\xi,}\\
 \vspace{0.1cm}
 \displaystyle{K_{0,2}(\sigma)(x)= \int_{\R^3} e^{ix\cdot \xi +i \frac{\sigma}{F} b_F(\xi)} \chi(\frac{|\xi_h|}{k_0}) \big(1-\chi(\frac{|\xi_3|}{k_0})\big) \varphi_1(|\xi|) d\xi,}\\
 \displaystyle{K_{0,3}(\sigma)(x)=\int_{\R^3} e^{ix\cdot \xi +i \frac{\sigma}{F} b_F(\xi)} \big(1-\chi(\frac{|\xi_3|}{k_0})\big) \big(1-\chi(\frac{|\xi_h|}{k_0})\big) \varphi_1(|\xi|) d\xi.}
\end{cases}
$$
Thanks to the study of $D^2 b_F$, we can apply Theorem \ref{ThLitt} (and more precisely the formulation from \eqref{ThLitt2}): there exists $C_F>0$ such that for all $j\in \{1,2,3\}$:
$$
\|K_{0,j}(\sigma)\|_{L^\infty} \leq C_F \min(1, |\sigma|^{-\frac{j}2}),
$$
which leads to:
$$
 \|K_0(\sigma)\|_{L^\infty} \leq C_F \sum_{j=1}^3 \min(1, |\sigma|^{-\frac{j}2}) \leq C_F \min(1, |\sigma|^{-\frac12}).
$$
Gathering this with \eqref{Chgtvar} and \eqref{EstimLjL1} leads to
$$
\|L_j(\sigma)g\|_{L^\infty}\leq C_F 2^{3j}\min(1, |\sigma|^{-\frac12})\|g\|_{L^1},
$$
In terms of norms, together with \eqref{EstimLjL2}, we can write for all $\theta \in[0,1]$:
$$
\begin{cases}
 \|L_j(\sigma)\|_{L^2 \rightarrow L^2} \leq C,\\
 \|L_j(\sigma)\|_{L^1 \rightarrow L^\infty} \leq C_F 2^{3j}|\sigma|^{-\frac{\theta}2},
\end{cases}
$$
and thanks to the Riesz-Thorin theorem, we end up, for all $q\in[2,\infty]$ and $\theta \in[0,1]$, with:
\begin{equation}
 \|L_j(\sigma)g\|_{L^q} \leq (C_F)^{1-\frac2{q}} \frac{2^{3j(1-\frac2{q})}}{|\sigma|^{\frac{\theta}2 (1-\frac2{q})}} \|g\|_{L^{\bar{q}}},
\end{equation}
and together with \eqref{Estimphi}, going back to \eqref{estimTT1}, we get:
$$
 \|\ddj f\|_{L^p L^r} \leq C_F^{\frac12-\frac1{r}} 2^{\frac{3j}2 (1-\frac2{r})} \ee^{\frac{\theta}4 (1-\frac2{r})} \|\ddj f_0\|_{L^2} \sup_{\psi \in \cB} \left[\int_0^\infty \int_0^\infty \frac{h(t)h(t')}{|t-t'|^{\frac{\theta}2 (1-\frac2{r})}}  dtdt'\right]^{\frac12},
$$
with $h(t)=e^{-\frac{\nu}4 t 2^{2j}} \|\psi(t,.)\|_{L^{\bar{r}}}$. Next, using the Hardy-Littlewood-Sobolev estimates as we did in \cite{FCPAA} (Proposition 49 with $\aa=\frac{\theta}{2}(1-\frac{2}{r})$, we also refer to \cite{HL, So, Lieb}), we get that
\begin{multline}
 \int_0^\infty \int_0^\infty \frac{h(t)h(t')}{|t-t'|^{\frac{\theta}{2}(1-\frac{2}{r})}} dt dt' \leq C\|h\|_{L^m}^2 \leq C\left( \|\mbox{exp}(-\frac{\nu}4  2^{2j} \cdot)\|_{L^k} \|\psi\|_{L^{\bar{p}}L^{\bar{r}}} \right)^2\\
 \leq C\left(\frac{1}{\nu^{\frac1{k}}} \left[\frac4{k}\right]^{\frac{1}{k}} 2^{-\frac{2j}{k}} \|\psi\|_{L^{\bar{p}}L^{\bar{r}}} \right)^2,
\end{multline}
for $m,k\in [1,\infty]$ chosen so that $\frac1{m}=1-\frac{\theta}4 (1-\frac2{r})$ and $\frac1{k}+\frac1{\bar{p}}=\frac1{m}$, that is:
$$
\frac{1}{k}=\frac{1}{p}-\frac{\theta}{4}(1-\frac{2}{r}).
$$
\begin{rem}
 \sl{As we want $k\geq 1$ we get the condition $p \leq \frac{4}{\theta (1-\frac{2}{r})}$.}
\end{rem}
We end up with:
$$
 \|\ddj f\|_{L^p L^r} \leq \frac{C_F^{\frac12-\frac1{r}}}{\nu^{\frac1{p}-\frac{\theta}4 (1-\frac{2}{r})}} \left[4\Big(\frac1{p}-\frac{\theta}4 (1-\frac{2}{r})\Big)\right]^{\frac1{p}-\frac{\theta}4 (1-\frac2{r})} \ee^{\frac{\theta}4 (1-\frac{2}{r})} 2^{j(\frac32-\frac{3}{r} -\frac{2}{p} +\frac{\theta}2 (1-\frac{2}{r}))} \|\ddj f_0\|_{L^2}
$$
which leads to the result in the homogeneous case. The inhomogeneous case (i.-e. when $\Fe\neq 0$) easily follows applying the very same steps to the Duhamel formula. $\blacksquare$

\section{Proof of Theorem \ref{Th3}}

As a consequence of Theorem \ref{Th1}, we immediately get the following estimates for the external force term (see Theorem 14 and Proposition 25 from \cite{FCPAA} for more details)
\begin{prop} \label{estimGlb}
 \sl{There exists a constant $C_F>0$ such that for all $\delta\in]0, \frac12]$, $\tUoqg\in H^{\frac12+\delta}(\R^3)$ (with $\max(1, \|\tUoqg\|_{H^{\frac12+\delta}(\R^3)})\leq \Co)$) and $s\in[0,\frac12+\delta]$ ($\Cdn$ being the same as in Theorem \ref{Th1}),
 \begin{equation}
\int_0^\infty \|G^b(\tau)\|_{\Hs } d\tau \leq \frac{C_F}{\nu} \Cdn \C_0^{2+\frac{1}{\delta}},
 \end{equation}
}
\end{prop}
Let us recall that (see \cite{FC2, FCPAA} for details) $\We$ is globally defined and that there exists a constant $\Do$ (only depending on $F,\nu, \delta, \Co$) such that for any $s\in[\frac12, \frac12 +\delta]$,
\begin{equation}
 \|\We\|_{\dot{E}^s}^2 \leq \Do \left(\|\Uoosc\|_{\Hs }^2 + 1\right),
\label{estimWe1}
 \end{equation}
Let us also recall that we already know the existence of a local strong solution $\Ue$ whose lifespan will be denoted as $T_\ee^*$. As $\tUqg$ and $\We$ exist globally, $\de$ is well defined in $\dot{E}_T^{\frac12} \cap \dot{E}_T^{\frac12 +\delta}$ for all  $T<T_\ee^*$ and we can perform the innerproduct in $\dot{H}^s$ of System \eqref{GE} with $\de$.

In order to bound $(F_4|\de)_{\dot{H}^s}$ to $(F_8|\de)_{\dot{H}^s}$, as we wish to get the biggest possible initial data, we had to use in \cite{FCPAA} a different approach from what we did in \cite{FC2, FC3} (in order to resort as little as possible to \eqref{estimWe1} allowing the largest possible initial data) involving the use of the following nonlocal operator: for $s\in]0,1[$, we define
\begin{multline}
|D|^s f(x)=C_s PV \int_{\R^3} \frac{f(x)-f(y)}{|x-y|^{3+s}} dy =C_s PV \int_{\R^3} \frac{f(x)-f(x-y)}{|y|^{3+s}} dy\\
= C_s \underset{a\rightarrow 0}{\mbox{lim}} \int_{|x-y|>a} \frac{f(x)-f(y)}{|x-y|^{3+s}} dy =C_s \underset{a\rightarrow 0}{\mbox{lim}} \int_{|y|>a} \frac{f(x)-f(x-y)}{|y|^{3+s}} dy
\end{multline}
Under appropriate assumptions the "PV" can be removed and we refer to \cite{FCpochesLp, FCPAA} for the following result:
\begin{prop}
\sl{For any $s\in]0,1[$ and any smooth functions $f,g$ we can write:
$$
|D|^s(fg)= (|D|^sf)g +f|D|^sg +M_s(f,g),
$$
where the bilinear operator $M_s$ is defined for all $x\in \R^3$ as:
\begin{equation}
M_s(f,g)(x) =\int_{\R^3} \frac{\big(f(x)-f(x-y)\big) \big(g(x)-g(x-y)\big)}{|y|^{3+s}} dy.
\label{defM}
\end{equation}
Moreover there exists a constant $C_s$ such that for all $f,g$ and all $r,r_1, r_2,q_1,q_2\in[1,\infty]$ and $s_1,s_2>0$ satisfying:
$$
\displaystyle{\frac{1}{r}= \frac{1}{r_1} +\frac{1}{r_2}, \quad 1= \frac{1}{q_1} +\frac{1}{q_2}, \quad s_1+s_2=s},
$$
then we have
\begin{equation}
\|M_s(f,g)\|_{L^r} \leq C_s \|f\|_{\displaystyle{\dot{B}_{r_1,q_1}^{s_1}}}\|g\|_{\displaystyle{\dot{B}_{r_2,q_2}^{s_2}}}.
\end{equation}
}
\label{propM2}
\end{prop}
As the computations are close for Theorems \ref{Th2} and \ref{Th3} we will first present the computations in some general $\dot{H}^s$ and then switch to the specificities of the case $s=\frac12$. For the sake of conciseness we will focus on the changes, which will only occur in the estimates for $F_8$ (we refer to Section 2.2.1 from \cite{FCPAA} for details about the other terms from \eqref{f1f2}). Let us recall that in \cite{FCPAA}, the best we could write was:
\begin{multline}
 \Big|(F_8\big|\de)_{\Hs}\Big| \leq \||D|^{s+\aC}(\We \otimes \We)\|_{L^\frac{6}{3+2\aC}} \||D|^{s-\aC} \n \de\|_{L^\frac{6}{3-2\aC}}\\
 \leq \left(2\||D|^{s+\aC}\We\|_{L^2} \|\We\|_{L^\frac{3}{\aC}} +\|\We\|_{\dot{H}^{s+\aC-\bd}}\|\We\|_{\dot{B}_{\frac{3}{\aC},2}^{\bd}} \right) \cdot \|\de\|_{\dot{H}^{s+1}}.
 \label{En8}
\end{multline}
Indeed the main constraints were imposed by the fact that we looked for the biggest possible bound for the initial oscillating part ($\ee^{-\gamma}$ with $\gamma<\frac{\delta}2$), and the fact that, using Proposition \ref{propM2} to estimate $M_{s+\aa_3}(\We, \We)$, we need:
$$
s_1, s_2>0 \quad \mbox{and} \quad \frac{1}{r_1} +\frac{1}{r_2}=\frac12+\frac{\aC}3.
$$
As in Proposition 47 from \cite{FCPAA} the integration index in the Besov spaces had to be strictly greater than 4, it was not possible to hope for two uses of the Strichartz estimates, as at least one index (among $r_1,r_2$) would then be stricly less than 4. We did what we tried at most to avoid: take $r_1=2$ and use the energy estimates for $\We$ together with Strichartz estimates for the other term ($r_2$ was then large as $\aC>0$ is small).
\\

For the first term in \eqref{En8} (as explained in (2.45) in \cite{FCPAA}), taking $(d,p,r,q,\theta)=(0,2,\frac{3}{\aC},2, \frac{6(\delta+\aC)}{3-4\aC})$ in the Strichartz estimates, for any $\aC\in ]0,\delta]$,
$$
\|\We\|_{L^\infty \dot{H}^{s+\aC}} \|\We\|_{L^2 L^\frac{3}{\aC}} \leq C_{F, \nu,\delta, \aC} \ee^{\frac{\delta+\aC}2} \Do (\|\Uoosc\|_{\dot{H}^{s+\aC}}^2+1)^\frac12 (\|\Uoosc\|_{\dot{H}^{\frac12+\delta}}+ 1).
$$
As we wished for the maximal possible power of $\ee$ our best option was to take $\aC=\delta+\frac12-s$ and then to estimate correctly the second term from \eqref{En8}, we needed that $\|\Uoosc\|_{\dot{H}^{\frac12}\cap \dot{H}^{\frac12+\delta}} \leq \Co \ee^{-\gamma}$ (in fact a bound for the $\dot{H}^{\frac12+c\delta}\cap \dot{H}^{\frac12+\delta}$ for some $c<1$ close to 1 was sufficient). We eventually needed $\delta< 3/26$.

With our new Strichartz estimates, the integration index is now allowed to be in $[2, \infty]$ so choosing in Proposition \ref{propM2}:
$$
r_1=r_2=\frac{12}{3+2\aC}, \quad q_1=q_2=2, \quad \mbox{and } s_1=s_2=\frac{s+\aC}2, 
$$
we can write ($r_1<4$), using that for all $p\in[2,\infty[$, $\dot{B}_{p,2}^0 \hookrightarrow L^p$ (see \cite{Dbook}):
\begin{multline}
 \Big|(F_8\big|\de)_{\dot{H}^s}\Big| \leq \Big(2\||D|^{s+\aC}\We\|_{L^\frac{12}{3+2\aC}} \|\We\|_{L^\frac{12}{3+2\aC}} +\|\We\|_{\dot{B}_{\frac{12}{3+2\aC}, 2}^{\frac{s+\aC}2}}^2 \Big) \cdot \|\de\|_{\dot{H}^{s+1}}\\
 \leq \Big(2\||D|^{s+\aC}\We\|_{\dot{B}_{\frac{12}{3+2\aC},2}^0} \|\We\|_{\dot{B}_{\frac{12}{3+2\aC},2}^0} +\|\We\|_{\dot{B}_{\frac{12}{3+2\aC}, 2}^{\frac{s+\aC}2}}^2 \Big) \cdot \|\de\|_{\dot{H}^{s+1}}\\
 \leq \frac{\nu}{16} \|\de\|_{\dot{H}^{s+1}}^2 +\frac{C}{\nu} \left(\||D|^{s+\aC}\We\|_{\dot{B}_{\frac{12}{3+2\aC},2}^0}^2 \|\We\|_{\dot{B}_{\frac{12}{3+2\aC},2}^0}^2 +\|\We\|_{\dot{B}_{\frac{12}{3+2\aC}, 2}^{\frac{s+\aC}2}}^4 \right),
 \label{En8bis}
\end{multline}
so that the last terms now become much nicer as we now do not use energy estimates for $\We$ \emph{at all}. Next, repeating the bootstrap arguments, if $\ee>0$ is so small that $\|\delta_\ee(0)\|_{\dot{H}\frac12} \leq \frac{\nu}{8C}$ (see (2.37) and (2.38) from \cite{FCPAA}) and if we define 
\begin{equation}
 T_\ee \overset{def}{=} \sup \{t\in[0,T_\ee^*[, \quad \forall t'\leq t, \|\de(t')\|_{\dot{H}^\frac12} \leq \frac{\nu}{4C}\},
\end{equation}
then we have, for all $t<T_\ee$:

\begin{multline}
 \|\de(t)\|_{\dot{H}^s}^2 +\frac{\nu}{2} \int_0^t \|\n \de(\tau)\|_{\dot{H}^s}^2 d\tau \leq \Bigg[ \|\de(0)\|_{\dot{H}^s}^2 +\frac{C}{\nu} \Bigg(\|\tUqg\|_{L^\infty \dot{H}^s}^{2(1-\aA)} \|\tUqg\|_{L^2 \dot{H}^{s+1}}^{2\aA} \|\We\|_{L^\frac{2}{1-\aA} L^\frac{3}{\aA}}^2\\
 +\|\tUqg\|_{L^\infty \dot{H}^\frac12}^2 \||D|^{s+\aA}\We\|_{L^2 L^\frac{6}{1+2\aA}}^2 +\|\tUqg\|_{L^\infty \dot{H}^s}^2 \|\We\|_{L^2 \dot{B}_{\frac{3}{\aA},2}^{\aA}}^2\\
 +\frac{\nu}{4C}\||D|^{s+\aB}\We\|_{L^2 L^\frac{6}{1+2\aB}}^2 +\||D|^{s+\aC}\We\|_{L^{p_1}\dot{B}_{\frac{12}{3+2\aC},2}^0}^2 \|\We\|_{L^{p_2}\dot{B}_{\frac{12}{3+2\aC},2}^0}^2 +\|\We\|_{L^4 \dot{B}_{\frac{12}{3+2\aC}, 2}^{\frac{s+\aC}2}}^4 \Bigg) \Bigg]\\
\times \exp \frac{C}{\nu}\Bigg\{ \|\n \tUqg\|_{L^2 \dot{H}^{\frac12}}^2 (1+\frac{1}{\nu^2}\|\tUqg\|_{L^\infty \dot{H}^{\frac12}}^2) +\frac{1}{\nu^{\frac{2\aB}{1-\aB}}}\|\We\|_{L^\frac{2}{1-\aB} L^\frac{3}{\aB}}^\frac{2}{1-\aB} +\|\We\|_{L^2 \dot{B}_{\frac{3}{\aB},2}^{\aB}}^2 \Bigg\},
\label{estimCas1}
\end{multline}
with $p_1,p_2$ to be fixed later and satisfying
\begin{equation}
 \frac1{p_1}+\frac1{p_2}=\frac12.
 \label{Condp1p2}
\end{equation}
Now, forcing $\sigma=\frac12 + \delta$ in Proposition \ref{Estimdispnu} (we refer to \cite{FCPAA} for details) and taking successively $(d,p,r,q, \theta)$ in
$$
\left\{(s+\aa,2,\frac{6}{1+2\aa},2,\frac{3}{1-\alpha}(\delta+\frac12-s)), \quad (\aa,2,\frac{3}{\aa},2,\frac{6\delta}{3-2\aa}), \quad (0,\frac{2}{1-\aa},\frac{3}{\aa},2,\frac{6\delta}{3-2\aa})\right\},
$$
we get:
\begin{equation}
\begin{cases}
\vspace{1mm}
 \||D|^{s+\aa} \We\|_{L_t^2 L^\frac{6}{1+2\aa}} \leq \Do \ee^{\frac12(\delta+\frac12-s)} (\|\Uoosc\|_{\dot{H}^{\frac12+\delta}}+ 1)\\
 \|\We\|_{L_t^2 \dot{B}_{\frac{3}{\aa},2}^\aa} +\|\We\|_{L_t^\frac{2}{1-\aa} L^\frac{3}{\aa}} \leq \Do \ee^{\frac{\delta}2} (\|\Uoosc\|_{\dot{H}^{\frac12+\delta}}+ 1).
\end{cases}
\label{StriA}
\end{equation}
\begin{rem}
 \sl{To bound the previous $L_T^a \dot{B}_{b,c}^s$-type norms, we use Proposition \ref{Propermut} to take advantage of the Strichartz estimates from Proposition \ref{Estimdispnu} which involve $\tilde{L}_T^a \dot{B}_{b,c}^s$-type norms. For example in the first estimate, as $p,q\geq 2$, we can write that $\|\We\|_{L^p L^q} \leq \|\We\|_{L^p \dot{B}_{q,2}^0} \leq \|\We\|_{\tilde{L}^p \dot{B}_{q,2}^0}$.}
\end{rem}

In each case, asking $\theta\in]0,1]$ and $p\leq \frac4{\theta(1-2/r)}$ leads to the conditions (for $\aA=\aB=\aa>0$ small):
\begin{equation}
\begin{cases}
\vspace{1mm}
\delta\leq 1-\aa,\\
\vspace{1mm}
 \delta+\frac12-s \leq \frac13-\frac{\aa}3,\\
 \delta \leq \frac12-\frac{\aa}3,
 \label{Precondelta2}
\end{cases}
\end{equation}
and as the second one has to be satisfied for any $s\in[\frac12, \frac12+2\eta_0 \delta]$, the conditions are verified when
\begin{equation}
 \delta \leq \min(\frac12-\frac{\aa}3, \frac13-\frac{\aa}3)=\frac13-\frac{\aa}3.
 \label{Condelta2}
\end{equation}
Let us continue with the last terms coming from the estimate of $F_8$: choosing $(d,p,r,q,\theta)=(\frac{s+\aC}2, 4, \frac{12}{3+2\aC}, 2,\frac{12\delta+6(\frac12-s)}{3-2\aC})$, we obtain:
\begin{equation}
 \|\We\|_{L^4 \dot{B}_{\frac{12}{3+2\aC}, 2}^{\frac{s+\aC}2}} \leq \Do \ee^{\frac12\left(\delta+\frac12(\frac12-s)\right)} (\|\Uoosc\|_{\dot{H}^{\frac12+\delta}}+ 1),
\label{StriB}
\end{equation}
and $\theta \leq 1$ (for any $s\in[\frac12, \frac12+2\eta_0 \delta]$) if and only if
\begin{equation}
\delta+\frac12 (\frac12-s)\leq \frac14-\frac{\aC}6,
\label{Precondelta3}
\end{equation}
which is true when
\begin{equation}
 \delta \leq \frac14-\frac{\aC}6.
 \label{Condelta3}
\end{equation}
The condition on p (in Proposition \ref{Estimdispnu}) leads to $\delta+\frac12 (\frac12-s)\leq \frac12$ which is obviously satisfied.
\begin{rem}
 \sl{We emphasize \eqref{StriB} is way better than (2.48) from \cite{FCPAA} in the case $s=\frac12$ (we will give details in what follows) but in the case $s\in]\frac12, \frac12 +2\eta_0 \delta]$ it does not allow anything better than $\gamma<\frac{\delta}2$.}
\end{rem}
Finally, let us precise $p_1$ and $p_2$ to estimate the last term: using Proposition \ref{Estimdispnu} for $(d,p,r,q)$ in
$$
\left\{(s+\aC, p_1, \frac{12}{3+2\aC}, 2), \quad (0, p_2, \frac{12}{3+2\aC}, 2)\right\},
$$
with (for some $\theta_{1,2}$ to be fixed later):
\begin{equation}
 \begin{cases}
 \sigma_1= s+\frac34 +\frac{\aC}2-\frac2{p_1}+\frac{\theta_1}2 (\frac12- \frac{\aC}3),\\
 \sigma_2= \frac34 -\frac{\aC}2-\frac2{p_2}+\frac{\theta_2}2 (\frac12- \frac{\aC}3).
\end{cases}
\label{Sigmas}
\end{equation}
If we want $\sigma_1=\sigma_2=\frac12 +\delta$, we can ask that:
$$
s+\frac{\aC}2-\frac2{p_1}=-\frac2{p_2}-\frac{\aC}2,
$$
which, combined with \eqref{Condp1p2} and \eqref{Sigmas}, leads to:
$$
(p_1,p_2)=(\frac4{1+s+\aC}, \frac4{1-s-\aC}), \quad \mbox{and } \theta_1=\theta_2=\frac{12\delta+6(\frac12-s)}{3-2\aC},
$$
and we obtain:
\begin{equation}
 \||D|^{s+\aC}\We\|_{L^\frac4{1+s+\aC} \dot{B}_{\frac{12}{3+2\aC}, 2}^0} +\|\We\|_{L^\frac4{1-s-\aC} \dot{B}_{\frac{12}{3+2\aC}, 2}^0}
 \leq \Do \ee^{\frac12(\delta+\frac12(\frac12-s))} (\|\Uoosc\|_{\dot{H}^{\frac12+\delta}}+ 1),
 \label{StriC}
\end{equation}
and asking $\theta_i\in]0,1]$ and $p_i\leq \frac4{\theta_i(1-2/r_i)}$ leads to:
\begin{equation}
\begin{cases}
\vspace{1mm}
\delta\leq \frac14-\frac{\aC}2,\\
 \delta+\frac12 (\frac12-s)\leq \frac14-\frac{\aC}6,
 \label{Precondelta4}
\end{cases}
\end{equation}
which is satisfied for all $s\in[\frac12, \frac12+2\eta_0 \delta]$ when
\begin{equation}
 \delta \leq \frac14-\frac{\aC}2.
 \label{Condelta4}
\end{equation}
Contrary to \cite{FCPAA} where we needed to impose $\aC=\delta+\frac12-s$, nothing prevents us from choosing $\aA=\aB=\aC=\aa\in]0,\frac12[$ small and combining \eqref{Condelta2}, \eqref{Condelta3} and \eqref{Condelta4}, we only need
\begin{equation}
\delta\leq \frac14-\frac{\aC}2.
 \label{Condelta}
\end{equation}
So that the only condition we impose will be $\delta<\frac14$ and we will choose $\aa=\frac12-2\delta$ (which garantees $s-\aa \in]0,1[$). Combining \eqref{estimCas1} with \eqref{StriA}, \eqref{StriB} and \eqref{StriC}, we end-up for all $t\leq T_\ee$ with
 \begin{multline}
 \|\de(t)\|_{\dot{H}^s}^2 +\frac{\nu}{2} \int_0^t \|\n \de(\tau)\|_{\dot{H}^s}^2 d\tau \leq\\
 \Bigg[ \|\de(0)\|_{\dot{H}^s}^2+\Do\Big( (\ee^\delta +\ee^{\delta+\frac12-s}) (\|\Uoosc\|_{\dot{H}^{\frac12+\delta}}+ 1)^2 +\ee^{2\delta+\frac12-s}(\|\Uoosc\|_{\dot{H}^{\frac12+\delta}}+ 1)^4\Big)\Bigg]\\
 \times \exp \Bigg\{\Do\left(1+\ee^\delta (\|\Uoosc\|_{\dot{H}^{\frac12+\delta}}+ 1)^2 +\Big(\ee^\delta (\|\Uoosc\|_{\dot{H}^{\frac12+\delta}}+1)^2\Big)^\frac1{1-\aa}\right)\Bigg\}.
\label{estimCas1s}
\end{multline}
Let us now switch to the case $s=\frac12$ to complete the proof of Theorem \ref{Th3}. The previous estimates turn into: 
 \begin{multline}
 \|\de(t)\|_{\dot{H}^\frac12}^2 +\frac{\nu}{2} \int_0^t \|\n \de(\tau)\|_{\dot{H}^\frac12}^2 d\tau\\
 \leq \Bigg[ \|\de(0)\|_{\dot{H}^\frac12}^2  +\Do\Big( \ee^\delta (\|\Uoosc\|_{\dot{H}^{\frac12+\delta}}+ 1)^2 +\ee^{2\delta}(\|\Uoosc\|_{\dot{H}^{\frac12+\delta}}+ 1)^4\Big)\Bigg]\\
 \times \exp \Bigg\{\Do\left(1+\ee^\delta (\|\Uoosc\|_{\dot{H}^{\frac12+\delta}}+ 1)^2 +\Big(\ee^\delta (\|\Uoosc\|_{\dot{H}^{\frac12+\delta}}^2+1)^2\Big)^\frac1{1-\aa}\right)\Bigg\}\\
 \leq \Do \left[\ee^{2\aa_0} +K_\ee +K_\ee^2\right] \exp \Big(\Do\Big(1+ K_\ee +K_\ee^\frac1{1-\aa} \Big)\Big)\\
 \leq \Do (\ee^{2\aa_0} +K_\ee) \exp \Big(\Do(1+ K_\ee)\Big),
\label{estimCas1b}
\end{multline}
where $K_\ee=\ee^\delta \|\Uoosc\|_{\dot{H}^{\frac12+\delta}}^2$. We recall that $\Uoosc$ is large but $K_\ee\leq m_0^2$ is small thanks to the assumptions so \eqref{estimCas1b} turns into:
\begin{equation}
 \|\de(t)\|_{\dot{H}^\frac12}^2 +\frac{\nu}{2} \int_0^t \|\n \de(\tau)\|_{\dot{H}^\frac12}^2 d\tau \leq \Do (\ee^{2\aa_0} +m_0^2) \exp \Big(\Do(1+ m_0^2)\Big),
\end{equation}
so if we ask $m_0\leq 1$ then \eqref{estimCas1b} is bounded by $\Do e^{2\Do}(\ee^{2\aa_0} +m_0^2)$, and when $\ee$ and $m_0$ are so small that
$$
\ee^{2\aa_0} +m_0^2\leq \Do^{-1}e^{-2\Do} \left(\frac{\nu}{8C}\right)^2
$$
then we can conclude the bootstrap argument and get $T_\ee=T_\ee^*=\infty$ (by the blow-up criterion \eqref{critereexpl}), which is the first result of Theorem \eqref{Th2}. If we ask that, in addition,
$$
\|\Uoosc\|_{\dot{H}^{\frac12+\delta}} \leq m(\ee) \ee^{-\frac{\delta}2},
$$
with $m(\ee)\rightarrow 0$ as $\ee$ goes to zero, then when $\ee$ is so small that $m(\ee)\leq 1$, \eqref{estimCas1b} is bounded by $\Do e^{2\Do}(\ee^{2\aa_0} +m(\ee)^2)$ and similar arguments allow us to prove global well-posedness and provide the following convergence rate:
$$
\|\de(t)\|_{\dot{H}^\frac12}^2 +\frac{\nu}{2} \int_0^t \|\n \de(\tau)\|_{\dot{H}^\frac12}^2 d\tau \leq \Do e^{2\Do} \left(\min(\ee^{\aa_0}, m(\ee))\right)^2,
$$
which concludes the proof of Theorem \ref{Th3}. $\blacksquare$

\section{Proof of Theorem \ref{Th2}}

To obtain the rest of Theorem \ref{Th2}, let us go back to \eqref{estimCas1s}, which turns into (taking $\aA=\aB=\aC>0$): for all $s\in [\frac12-\eta \delta, \frac12+\eta \delta]$ and $t<T_\ee$,
\begin{equation}
\|\de(t)\|_{\dot{H}^s}^2 +\frac{\nu}{2} \int_0^t \|\n \de(\tau)\|_{\dot{H}^s}^2 d\tau \leq \Do \left[\ee^{2\aa_0} +\ee^{\delta-2\gamma} +\ee^{\delta-2\gamma-(s-\frac12)}\right] e^{\Do(1+\ee^{\delta-2\gamma})}.
 \label{estimCas1sb}
\end{equation}
Following the lines of the end of Section 2.2 from \cite{FCPAA} leads to the conclusion of the boostrap method and provides the first point. But prior to that, as the previous estimates are for $s\in [\frac12-\eta \delta, \frac12+\eta \delta]$, the conditions coming from the use of our Strichartz estimates, namely \eqref{Precondelta2}, \eqref{Precondelta3} and \eqref{Precondelta4} (imposed by the fact that all the $\theta$ parameters defined in the use of Proposition \ref{Estimdispnu} lie in ]0,1] and satisfy $p\leq \frac4{\theta(1-\frac2{r})}$) will have to be fulfilled for this extended range for $s$, and we need that (the first condition implies that $0<\aC<\min(s,1-s)$ which is equivalent to $s\pm \aC \in ]0,1[$):
 $$
 \begin{cases}
\vspace{1mm}
0<\aC <\frac12-\eta \delta,\\
\vspace{1mm}
 (1+\eta)\delta \leq \frac13-\frac{\aC}3,\\
 \vspace{1mm}
 (1+\frac{\eta}2)\delta \leq \frac14-\frac{\aC}6,\\
 \delta \leq \frac14-\frac{\aC}2.
\end{cases}
 $$
Let us first ask that $\delta<\frac14$ and put $\beta_0=\frac12-2\delta$ (so that $\delta=\frac14-\frac{\beta_0}2$), then define $\gamma$ according to \eqref{defdeltaeta}. The last condition asks $0<\aC\leq\beta_0$, so choosing $\aC=\frac{\beta_0}{K}$ for some large $K\geq 1$, if $K>4\beta_0$ then the first condition is implied by the second one ($\delta< \frac14$ and $\eta\leq 1$) and therefore the condition on $\eta\in ]0,2\eta_0[$ is implied by:
$$
 0< \eta\leq \min\Big(\frac{\frac13+\beta_0(2-\frac4{3K})}{1-2\beta_0}, \frac{2\beta_0(2-\frac2{3K})}{1-2\beta_0}\Big),
$$
Returning to the original parameters, the condition is fulfilled when
$$
 0< \eta\leq \min\Big(1-\frac{2\gamma}{\delta}, \frac1{12\delta}+(\frac1{4\delta}-1)(1-\frac2{3K}), (\frac1{4\delta}-1)(2-\frac2{3K})\Big),
$$
So, when $0<\eta <\min\big(1-2\frac{\gamma}{\delta}, \frac1{3\delta}-1, \frac1{2\delta}-2\big)$, choosing $K>\max(1, 2-8\delta)$ sufficiently large and taking $\aA=\aB=\aC=\frac2{K}(\frac14-\delta)$ fulfills the required conditions and allows us to consider \eqref{estimCas1sb} for all $t\leq T_\ee$ and all $s\in [\frac12-\eta \delta, \frac12+\eta \delta]$. When $\ee$ is so small that $\ee^{2 \eta_0 \delta} +(\ee^{2 \eta_0 \delta})^{\frac1{1-\aa}}\leq 1$ we get that:
\begin{multline}
 \|\de\|_{\dot{E}^s}^2 \leq \Do \left(\ee^{2\alpha_0} +\ee^{\delta-2\gamma} +\ee^{\delta-2\gamma-(s-\frac12)} +\ee^{2\delta-4\gamma-(s-\frac12)} \right) e^{\Do (1+\ee^{2 \eta_0 \delta} +(\ee^{2 \eta_0 \delta})^{\frac1{1-\aa}})} \\
\leq \Do e^{2\Do} \ee^{\min \big(2\alpha_0, \delta-2\gamma, \delta-2\gamma -(s-\frac12)\big)},
\end{multline}
and the rest of the proof of point 1 follows (as in \cite{FCPAA} or in the previous section, if $\ee$ is so small that the right-hand-side is less than $\nu/8C$ we obtain that $T_\ee=T_\ee^*=+\infty$).

Let us continue with the proof of point 3: as in \cite{FCPAA}, using Lemma \ref{majBs21} with $\aa=\beta=\eta \delta$, and the previous estimates for $s\in\{\frac12-\eta \delta, \frac12 +\eta \delta\}$, we deduce:
\begin{multline}
 \|\de\|_{L^2 L^\infty}\leq \|\de\|_{L^2 \dot{B}_{\infty,1}^0} \leq \|\de\|_{L^2 \dot{B}_{2,1}^\frac32}\leq C(\|\de\|_{L^2 \dot{H}^{\frac32-\eta \delta}} \|\de\|_{L^2 \dot{H}^{\frac32+\eta \delta}})^\frac12\\
 \leq \Bo \left(\ee^{\min(\aa_0, \eta_0\delta)} \ee^{\min(\aa_0, (2\eta_0-\eta)\frac{\delta}{2})}\right)^\frac12.
\label{estimdeLinf}
\end{multline}
Choosing $(d,p,r,q)=(0,2,\infty,1)$ and for all $\theta\in[0,1]$, from Proposition \ref{Estimdispnu}:
$$
 \|\We\|_{L^2 L^\infty} \leq \Co \ee^\frac{\theta}{4} \left(\|\Uoosc\|_{\dot{B}_{2,1}^{\frac12 +\frac{\theta}{2}}}+\int_0^\infty \|G^b(\tau)\|_{\dot{B}_{2,1}^{\frac12 +\frac{\theta}{2}}} d\tau\right).
$$
Using once more Lemma \ref{majBs21} with $(\aa,\beta)=(a \frac{\theta}{2}, b \frac{\theta}{2})$, and if $\theta=\frac{2\delta}{1+b}$ and $a=1-c(1+b)$(for some small $a,b>0$),
\begin{multline}
 \|\Uoosc\|_{\dot{B}_{2,1}^{\frac12 +\frac{\theta}{2}}}\leq \|\Uoosc\|_{\dot{H}^{\frac12+(1-a)\frac{\theta}2}}^\frac{b}{a+b} \|\Uoosc\|_{\dot{H}^{\frac12 +(1+b)\frac{\theta}{2}}}^\frac{a}{a+b}\\
 \leq \|\Uoosc\|_{\dot{H}^{\frac12+c \delta}}^\frac{b}{a+b} \|\Uoosc\|_{\dot{H}^{\frac12 +\delta}}^\frac{a}{a+b} \leq \Co \ee^{-\gamma}.
\end{multline}
Choosing $b=\frac{\eta'}{1-\eta'}$ (for $\eta'\in]0,1[$), we get
\begin{equation}
 \|\We\|_{L^2 L^\infty} \leq \Do \ee^{\frac{\delta}{2}(\frac{1}{1+b}-(1-2\eta_0))}=\Do \ee^{(2\eta_0-\eta') \frac{\delta}{2}},
 \label{estimweLinf}
\end{equation}
Asking that $a>0$ is equivalent to $\eta' <1-c$, so if we choose $\eta \in]0,\min\big(1-2\frac{\gamma}{\delta}, \frac1{4\delta}-1\big)[$ and $\eta'\in]0, \min(1-2\frac{\gamma}{\delta}, 1-c)[$, the conclusion follows from \eqref{estimdeLinf}, \eqref{estimweLinf} and the fact that we can take $\eta, \eta'>0$ as small as we want:
\begin{multline}
\|\Ue-\tUqg\|_{L^2 L^\infty}=\|\de +\We\|_{L^2 L^\infty}\\
\leq \Bo \ee^{\frac12\big(\min(\aa_0, \eta_0\delta) +\min(\aa_0, (2\eta_0-\eta)\frac{\delta}{2})\big)} +\Bo \ee^{(2\eta_0-\eta') \frac{\delta}{2}} \leq \Bo \ee^{\min(\aa_0, k\eta_0\delta)},
\end{multline}
for any $k<1$ as close to $1$ as we wish. $\blacksquare$
\\

Let us end with the proof of point 2: for any $\eta\in]0,2\eta_0[$, any $\eta'\in]0, \eta[$ and using the previous estimates for $s\in\{\frac12, \frac12 +\eta \delta\}$, we get:
\begin{multline}
 \big\||D|^{\eta' \delta}\de\big\|_{L^2 L^\infty}\leq \big\||D|^{\eta' \delta}\de\big\|_{L^2 \dot{B}_{2,1}^\frac32}\leq C \big\||D|^{\eta' \delta}\de\big\|_{L^2 \dot{H}^{\frac32-\eta' \delta}}^{1-\frac{\eta'}{\eta}} \big\||D|^{\eta' \delta}\de\big\|_{L^2 \dot{H}^{\frac32+(\eta-\eta') \delta}}^{\frac{\eta'}{\eta}}\\
 \leq C \|\de\|_{L^2 \dot{H}^{\frac32}}^{1-\frac{\eta'}{\eta}} \|\de\|_{L^2 \dot{H}^{\frac32+\eta \delta}}^{\frac{\eta'}{\eta}}\leq \Bo \ee^{\left((1-\frac{\eta'}{\eta})\min(\aa_0, \eta_0\delta) +\frac{\eta'}{\eta} \min(\aa_0, (\eta_0-\frac{\eta}2)\delta)\right)}.
\label{estimdeLinfb}
\end{multline}
Doing as in the previous step, but choosing $(d,p,r,q)=(\eta' \delta,2,\infty,1)$, we get that for all $\theta\in[0,1]$, from Proposition \ref{Estimdispnu}:
$$
 \big\||D|^{\eta' \delta}\We\big\|_{L^2 L^\infty} \leq \Co \ee^\frac{\theta}{4} \left(\|\Uoosc\|_{\dot{B}_{2,1}^{\frac12 +\eta' \delta +\frac{\theta}{2}}}+\int_0^\infty \|G^b(\tau)\|_{\dot{B}_{2,1}^{\frac12 +\eta' \delta +\frac{\theta}{2}}} d\tau\right).
$$
We then take (in Lemma \ref{majBs21}) $(\aa,\beta)=(a \frac{\theta}{2}, b \frac{\theta}{2})$ with $a,b$ such that
$$\begin{cases}
   \frac{\theta}2(1-a)+\eta'\delta=c\delta,\\
   \frac{\theta}2(1+b)+\eta'\delta=\delta.
  \end{cases}
$$
Such $(\aa,\beta)$ exist, we just need to take some $\eta"\in]0, \frac{1-c}{1-\eta'}[$, then define $b=\frac{\eta"}{1-\eta"}$, $a=1-(1+b)\frac{c-\eta'}{1-\eta'}$ and $\theta=2\delta (1-\eta')(1-\eta")$ ($c<1$ is close to 1, and the fact that $\eta'<c$ is implied by the first condition). With these choices, we end up with (the last estimate in \eqref{estimweLinfc} being valid if $0<\eta"\leq \frac{\eta-\eta'}{1-\eta'}$):
\begin{equation}
 \big\||D|^{\eta' \delta}\We\big\|_{L^2 L^\infty} \leq \Do \ee^{\frac{\delta}{2}(1-\eta')(1-\eta") -\gamma}=\Do \ee^{(\eta_0-\frac{\eta'}2-\frac{\eta"}2 (1-\eta')) \delta} \leq \Do \ee^{(\eta_0-\frac{\eta}2)\delta}.
 \label{estimweLinfc}
\end{equation}
Since
\begin{multline}
 (1-\frac{\eta'}{\eta})\min(\aa_0, \eta_0\delta) +\frac{\eta'}{\eta} \min(\aa_0, (\eta_0-\frac{\eta}2)\delta)\\
 \geq (1-\frac{\eta'}{\eta})(\eta_0-\frac{\eta}2)\delta) +\frac{\eta'}{\eta} \min(\aa_0, (\eta_0-\frac{\eta}2)\delta) =(\eta_0-\frac{\eta}2)\delta,
\end{multline}
we end-up with
$$
\big\||D|^{\eta'\delta}(\Ue-\tUqg)\big\|_{L^2 L^\infty} \leq \Bo \ee^{(\eta_0-\frac{\eta}2)\delta},
$$
which ends the proof. $\blacksquare$

\section{Appendix}

\subsection{Notations, Sobolev spaces and Littlewood-Paley decomposition}

We refer to the appendix of \cite{FCPAA} for general notations and properties of the Sobolev spaces as well as the Littlewood-Paley decomposition and the main properties we use. For a complete presentation, we refer to \cite{Dbook}.

Let us mention the following lemma whose proof is close to Lemma $5$ from \cite{FCestimLp} (see also Section 2.11 in \cite{Dbook}):
\begin{lem} \label{majBs21}
 \sl{For any $\aa, \beta>0$ there exists a constant $C_{\aa, \beta}>0$ such that for any $u\in \dot{H}^{s-\aa} \cap \dot{H}^{s+\beta}$, then $u\in\dot{B}_{2,1}^s$ and:
\begin{equation}
 \|u\|_{\dot{B}_{2,1}^s} \leq C_{\aa, \beta} \|u\|_{\dot{H}^{s-\aa}}^{\frac{\beta}{\aa + \beta}} \|u\|_{\dot{H}^{s+\beta}}^{\frac{\aa}{\aa + \beta}}.
\end{equation}
 }
\end{lem}
\begin{prop}
 \sl{\cite{Dbook} We have the following continuous injections:
$$
 \begin{cases}
\mbox{For any } p\geq 1, & \dot{B}_{p,1}^0 \hookrightarrow L^p,\\
\mbox{For any } p\in[2,\infty[, & \dot{B}_{p,2}^0 \hookrightarrow L^p,\\
\mbox{For any } p\in[1,2], & \dot{B}_{p,p}^0 \hookrightarrow L^p.
\end{cases}
$$
}
 \label{injectionLr}
\end{prop}
We end this section by the following proposition dedicated to Chemin-Lerner time-space Besov spaces: instead of considering $L^p \dot{B}_{q,r}^s$-type estimates, the integration in time is performed before the summation with respect to the frequency decomposition index
\begin{defi} \cite{Dbook}
 \sl{For $s,t\in \R$ and $a,b,c\in[1,\infty]$, we define the following norm
 $$
 \|u\|_{\tilde{L}_t^a \dot{B}_{b,c}^s}= \Big\| \left(2^{js}\|\ddj u\|_{L_t^a L^b}\right)_{j\in \Z}\Big\|_{l^c(\Z)}.
 $$
 The space $\tilde{L}_t^a \dot{B}_{b,c}^s$ is defined as the set of tempered distributions $u$ such that $\lim_{j \rightarrow -\infty} S_j u=0$ in $L^a([0,t],L^\infty(\R^d))$ and $\|u\|_{\tilde{L}_t^a, \dot{B}_{b,c}^s} <\infty$.
 }
 \label{deftilde}
\end{defi}
We refer once more to \cite{Dbook} (Section 2.6.3) for more details and will only recall that we have the following relations:
\begin{prop}
\sl{
For all $a,b,c\in [1,\infty]$ and $s\in \R$:
     $$
     \begin{cases}
    \mbox{if } a\leq c,& \forall u\in L_t^a \dot{B}_{b,c}^s, \quad \|u\|_{\tilde{L}_t^a \dot{B}_{b,c}^s} \leq \|u\|_{L_t^a \dot{B}_{b,c}^s}\\
    \mbox{if } a\geq c,& \forall u\in\tilde{L}_t^a \dot{B}_{b,c}^s, \quad \|u\|_{\tilde{L}_t^a \dot{B}_{b,c}^s} \geq \|u\|_{L_t^a \dot{B}_{b,c}^s}.
     \end{cases}
     $$
     \label{Propermut}
     }
\end{prop}

\subsection{A remark about the Strichartz estimates}
\label{RemLitt}
As outlined previously, had we not used the Littman theorem, we would end-up with the very close following result (whose proof only differs from the proof of Proposition 47 from \cite{FCPAA} in that we use the Riesz-Thorin theorem):
\begin{prop}
 \sl{There exists a constant $C_F>0$, such that for any $d\in \R$, $r>2$, $q\geq 1$, $\theta\in[0,\frac12[$ and $p\in[1, \frac{4}{\theta (1-\frac{2}{r})}]$, if $f$ solves \eqref{systdisp} for initial data $f_0$ and external force $\Fe$ both with zero divergence and potential vorticity, then ($c_0$ refers to the smaller constant appearing in the Littlewood-Paley decomposition, usually $c_0=\frac34$.)
 \begin{equation}
  \||D|^d f\|_{\tilde{L}_t^p\dot{B}_{r, q}^0} \leq \frac{C_F C_{p,\theta,r}}{\nu^{\frac{1}{p}-\frac{\theta}{4}(1-\frac{2}{r})}} \ee^{\frac{\theta}{4}(1-\frac{2}{r})} \times \left( \|f_0\|_{\dot{B}_{2, q}^\sigma} +\int_0^t  \|\Fe(\tau)\|_{\dot{B}_{2, q}^\sigma} d\tau \right),
 \end{equation}
 where
 $$
 \begin{cases}
  \sigma= d+\frac32-\frac{3}{r}-\frac{2}{p}+\frac{\theta}{2} (1-\frac{2}{r}),\\
  C_{p,\theta,r}=\left[\frac{2}{c_0^2}\Big(\frac{1}{p}-\frac{\theta}{4}(1-\frac{2}{r})\Big)\right]^{\frac{1}{p}-\frac{\theta}{4}(1-\frac{2}{r})} \displaystyle{\frac{2^{\frac12 (1-2\theta) (1-\frac{2}{r})}}{1-2^{-\frac12 (1-2\theta) (1-\frac{2}{r})}}}.
 \end{cases}
 $$
 }
\end{prop}
\begin{rem}
 \sl{\begin{enumerate}
      \item The differences for $\theta, r$ (in Proposition \ref{Estimdispnu}, $r\geq 2$ and $\theta \in[0,1]$) and the expression of $C_{p,\theta,r}$ only come from the fact that if we reproduce the proof of Proposition 47 from \cite{FCPAA} and begin with the vertical Littlewood-Paley truncation ($\ddk u=\varphi(2^{-j} D_3)u$):
$$
\|\ddj f\|_{L_t^p L_x^r}= \|\ddj f\|_{L^p L^r} \leq \Sum_{k=-\infty}^{j+1} \|\ddk \ddj f\|_{L^p L^r},
$$
we will have, in the end, to sum (for $k\leq j+1$) the term $2^{\frac{k-j}2 (1-2\theta) (1-\frac{2}{r})}$ which requires $(1-2\theta) (1-\frac{2}{r})>0$. Using this proposition instead of Proposition \ref{Estimdispnu} would also lead to the same Theorems \ref{Th2} and \ref{Th3} except that we would require $\delta<\frac18$ instead of $\delta< \frac14$.
\item We recall that in our case (shared with \cite{IMT} and also \cite{LT} in the case with only stratification) using the Littman theorem gives the same power $\left(\frac{\ee}{|t-t'|}\right)^\frac12$ as our first method in \cite{FCPAA} because one of the truncation (namely $K_{0,1}(\sigma)$) presents only one nonzero eigenvalue in the Hessian of the phase.
\item In the case of the rotating fluids (see Remark \ref{remVPRF} and \cite{KLT}) the Littman method allows to get $\frac{\ee}{|t-t'|}$ which improves the power of $\ee$ in the Strichartz estimates into $\frac{\theta}{2}(1-\frac{2}{r})$, but we pay the same "additional regularity" price as before as $\sigma= d+\frac32-\frac{3}{r}-\frac{2}{p}+\theta(1-\frac{2}{r})$.

For instance, for $(d,p,r,q)=(\aa, 2, \frac3{\aa},2)$, we get $\ee^{\frac{\theta}{4}(1-\frac{2\aa}3)}$ and $\sigma=\frac12+\frac{\theta}{2}(1-\frac{2\aa}3)$ in the first case, and to $\ee^{\frac{\theta}{2}(1-\frac{2\aa}3)}$ and $\sigma=\frac12+\theta(1-\frac{2\aa}3)$ if we use the Littman theorem. As the initial regularity is $\frac12+\delta$ for a small $\delta$, we reach the power $\ee^{\frac{\delta}2}$ in both cases taking successively $\theta$ in $\{\frac{2\delta}{1-\frac{2\alpha}3}, \frac{\delta}{1-\frac{2\alpha}3}\}$, which is possible as $\delta, \aa$ are small. So we eventually reach the same result and the use of the Littman theorem for our considerations does not give any improvement to the power of $\ee$.
\item In the case of the rotating fluids system, a visible difference between these two methods would be clearer when estimating the whole solution, as seen as a inhomogeneous solution of the linearized (dispersive) equation.
\end{enumerate}
}
 \label{plusdek}
\end{rem}

\textbf{Aknowledgements :} This work was supported by the ANR project INFAMIE, ANR-15-CE40-0011. We thank the anonymous referee for useful suggestions.

\end{document}